\documentclass[12pt]{amsart}
\usepackage[english]{babel}
\usepackage[applemac]{inputenc}
\usepackage[T1]{fontenc}

\usepackage{bbm}
\usepackage{float, verbatim}
\usepackage{amsmath}
\usepackage{amssymb}
\usepackage{amsthm}
\usepackage{amsfonts}
\usepackage{graphicx}
\usepackage{mathtools}
\usepackage{enumitem}
\usepackage[pagebackref]{hyperref}
\usepackage{color}
\usepackage{tikz}
\usepackage{xcolor}
\usepackage{pgfplots}
\usepackage{caption}
\usepackage{subcaption}
\usepackage{enumerate}
\usepackage{url}

\hypersetup{
	colorlinks=true,
	linkcolor=blue,
	filecolor=magenta,      
	urlcolor=cyan,
	citecolor= red,
}

\usepackage[final]{showlabels}

\newcommand{\ubox}{\overline{\dim}_{\mathrm{B}}}

\newcommand{\lbox}{\underline{\dim}_{\mathrm{B}}}
\newcommand{\boxd}{\dim_{\mathrm{B}}}

\newcommand{\Haus}{\dim_{\mathrm{H}}}

\newtheorem*{thm*}{Theorem}
\newtheorem*{conj*}{Conjecture}
\newtheorem*{ques*}{Question}
\newtheorem*{rem*}{Remark}
\newtheorem*{defn*}{Definition}
\newtheorem*{mainques*}{Main questions}

\newtheorem{thmx}{Theorem}

\newtheorem{thm}{Theorem}[section]

\newtheorem{cor}[thm]{Corollary}
\newtheorem{defn}[thm]{Definition}

\newtheorem{conj}[thm]{Conjecture}
\newtheorem{rem}[thm]{Remark}
\newtheorem{ques}[thm]{Question}

\newtheorem{exm}[thm]{Example}

\def\supp{\mathrm{supp}}

\def\cl{\mathrm{cl}}

\pgfplotsset{compat=1.16}

\begin{document}
\title{Missing digits points near manifolds}
%    author one information

\author{Han Yu}
\address{Han Yu, Mathematics Institute, Zeeman Building, University of Warwick, Coventry CV4 7AL, UK}
\curraddr{}
\email{han.yu.2@warwick.ac.uk}

%    \subjclass is required.


\maketitle

\begin{abstract}
We consider a problem concerning the distribution of points with missing digits coordinates that are close to non-degenerate analytic submanifolds. We show that large enough (to be specified in the paper) sets of points with missing digits coordinates distribute 'equally' around non-degenerate submanifolds. As a consequence, we show that intersecting those missing digits sets with non-degenerate submanifolds always achieve the optimal dimension reduction. On the other hand, we also prove that there is no lack of points with missing digits that are contained in non-degenerate submanifolds. Among the other results, 
\begin{itemize}
    \item[1]  we prove that the pinned distance sets of those missing digits sets contain non-trivial intervals regardless of where the pin is. 
\item[2] we prove that for each $\epsilon>0,$ for missing digits sets $K$ with large bases, simple digit sets (to be specified in the paper), and $\Haus K>3/4+\epsilon,$ the arithmetic product sets $K\cdot K$ contains non-trivial intervals.
\end{itemize}
\end{abstract}

\maketitle
%\tableofcontents
\allowdisplaybreaks
\section{Introduction}
We discuss a problem concerning the distribution of missing digits points around as well as on manifolds. Before we state the most general results, we list three special cases with some mixtures of number theory and geometric measure theory. 

In what follows, let $K_1$ be the set of points on $[0,1]$ whose base $10^{9000}$ expansions\footnote{Examples with much smaller bases can be found. However, we did not carry out the most optimal numerical computations.} contain only digits in $\{0,\dots,10^{8100}-1\}.$ Let $K_2$ be the set of points on $[0,1]$ whose base $10^{9000}$ expansions contain only digits in $\{0,\dots,10^{7000}-1\}.$ We see that $\Haus K_1=9/10$ and $\Haus K_2=7/9.$

\begin{thmx}\label{thm: special 1}
	Let $(t,t^2,t^3)_{t\in\mathbb{R}}$ be the Veronese curve in $\mathbb{R}^3.$ There is an integer $l\geq 0$ such that there are infinitely many $t>0$ such that the fractional parts of
	\[
	10^{9000l} t, 10^{9000l} t^2, 10^{9000l} t^3
	\]
	are contained in $K_1.$ Loosely speaking, there are many points on the Veronese curve whose coordinates have special expansions in base $10^{9000}.$ Moreover, the upper box dimension of the set consisting such numbers $t\in [0,1]$ is in $
	[1/30,7/10].$
\end{thmx}
\begin{rem*}
We expect that the (upper) box dimension of such $t$ should be exactly equal to $7/10.$ This is obtained via
\[
7/10=3*9/10-(3-1).
\]
Here $3*9/10$ is the dimension of the missing digits set $K_1\times K_1\times K_1$ in $\mathbb{R}^3$ considered in this theorem, $1$ is the dimension of the Veronese curve. For more discussions, see Theorem \ref{Main2}.
\end{rem*}
\begin{thmx}\label{thm: special 2}
	Let $K=K_2\times K_2\subset [0,1]^2$ be the twofold Cartesian product of $K_2.$ Then for each $x\in\mathbb{R}^2,$ the pinned distance set
	\[
	\Delta_x(K)=\{|x-y|:y\in K\}
	\]
	contains non-trivial interior and in particular, positive Lebesgue measure. Moreover, for each circle, $C\subset\mathbb{R}^2,$ we have
	\[
	\ubox C\cap K\leq \Haus K-(2-\dim C)=\frac{5}{9}.
	\]
	As a consequence, we can take $x=(0,0)$ and see that the set  \[K^2_2+K^2_2=\{x^2+y^2: x,y\in K_2\}\] contains intervals.
\end{thmx}
\begin{rem*}
Previously, it is possible to show that:

1. $\Delta_x(K)$ has positive Lebesgue measure for many choices of $x\in\mathbb{R}^2.$  Notice that the set $K$ has Hausdorff dimension larger than $1.5$. This statement follows from a classical result of Falconer. See \cite[Section 7]{Ma2}.

2. $\Delta_x(K)$ has full Hausdorff dimension for all $x\in\mathbb{R}^2.$ This follows by adapting the arguments in \cite{H14}.

This theorem pushes the result further, i.e. $\Delta_x(K)$ has a positive Lebesgue measure for all $x\in\mathbb{R}^2.$  This result holds in $\mathbb{R}^n,n\geq 3$ as well. Also, instead of a missing digits set in $\mathbb{R}^2,$ $K$ can be also chosen to be the Cartesian product of two missing digits sets with possibly different bases on $\mathbb{R}.$ See Section \ref{sec: application}.
\end{rem*}
\begin{thmx}\label{thm: special 3}
Consider the set $K_2.$ The arithmetic product set
\[
K_2\cdot K_2=\{xy:(x,y)\in K_2\times K_2\}
\]
has non-trivial interior.
\end{thmx}
\begin{rem*}
The sum set $K_2+K_2$ does not contain non-trivial intervals because it is contained in the set of numbers whose base $10^{9000}$ expansions contain only digits in $\{0,\dots,2\times 10^{7000}-2\}.$ Thus the nonlinearity of the multiplication map makes a significant difference.
\begin{comment}
By multiplying a large power of $10,$ we see that there is an integer $a$ such that
\[
\cup_{0\leq t_1,t_2\leq a} (t_1+K_2)(t_2+K_2)
\]
contains intervals. This union is a finite union, thus there is at least one choice $t_1,t_2$ such that
$
(t_1+K_2)(t_2+K_2)
$
already contains intervals.

In the theorem, it is likely that $l=0.$ In fact, this would follow from a non-rigorous statement in Remark \ref{Remark: unproved}(b). Formulating and proving precisely this statement would need a great amount of extra technical arguments. We decide to only prove this slightly weaker theorem for simplicity. We hope that this is a good trade-off. 
\end{comment}
\end{rem*}

The role of missing digits sets in this paper is by no means restrictive. It really just depends on Theorem \ref{thm: l1 bound}. In fact, it is possible to extend Theorem \ref{thm: l1 bound} to homogeneous self-similar measures with rational scaling ratios, finite rotation groups and rational translation vectors. However, in this paper, we anyway restrict the discussions to missing digits measures because of their intrinsic connections with number theory.

The roles of the specific missing digits sets $K_1,K_2$ with base $10^{9000}$ are also not restrictive. The precise conditions are on their Fourier $l^1$-dimensions, see Section \ref{sec: Fourier Dimension}. We need that $\dim_{l^1} K_1>8/9$ and $\dim_{l^1} K_2>3/4.$ These conditions are in some sense optimal. See Theorem \ref{thm: fail} and Remark \ref{rem: fail}.

For missing digits sets with simple digit sets (e.g. with consecutive digits), their Fourier $l^1$-dimensions are almost equal to their Hausdorff dimensions, see Theorem \ref{thm: l1 bound}.  Theorem \ref{thm: special 3} can be seen as a partial answer to the following folklore conjecture.
\begin{conj}\label{conj: Palis}
    Let $K$ be a missing digits set with $\Haus K>1/2,$ then $K\cdot K$ contains intervals.
\end{conj}
\begin{rem}\label{NewhouseRemark}
Let's see why this conjecture is plausible. By adapting the arguments in \cite{H14}, it can be shown that $\Haus K\cdot K=1$. Due to the nonlinearity of the multiplication map, it is expected that $K\cdot K$ has a positive measure. Finally, due to the results in \cite{MY}, it is expected that $K\cdot K$ contains intervals. 
\end{rem}
\begin{rem}
If $K$ has Newhouse thickness (see \cite[Definition 3.6]{KT}) $\tau(K)$ larger than $1$, then $K\cdot K$ contains intervals. In fact, with the help of Newhouse thickness of Cantor set (\cite{HKY},\cite{ZRZJ},\cite{KT},\cite{T}), it is possible to prove some of the topological results in Theorems \ref{thm: special 1}, \ref{thm: special 2}, \ref{thm: special 3}  for Cantor sets with thickness at least one. To see where $K_1,K_2$ are in the thickness story: $K_1$ at the beginning has Newhouse thicknesses $1/(10^{900}-1-10^{-900});$ $K_2$ has Newhouse thickness $1/(10^{2000}-1-10^{-2000}).$ Thus, in the sense of Newhouse thickness, they are very thin.
\end{rem}
\begin{rem}
By \cite{MS18}, it is possible to prove that for each missing digits set $K$, there is an integer $N$ so that the $N$-fold multiplication set $K\cdot K\cdot\dots\cdot K$ contains intervals. In fact, it is even possible to show that the corresponding measure $m_*((\mu_K)^N)$ is a smooth function where $m:(x_1,\dots,x_N)\to x_1\dots x_N$ is the multiplication map and $\mu_K$ is the natural missing digits measure supported on $K.$
\end{rem}
To motivate the general results in this paper, we first introduce three open problems (Questions \ref{ques0}, \ref{ques1}, \ref{ques2}). They are all linked together with a generalized version of a difficult conjecture of \'{E}. Borel. See \cite{Borel}. At the end of this paper, we will provide some more specific applications of our general results. See Sections \ref{sec: application}, \ref{sec: fractal measure}.

\subsection{Missing digits solutions for algebraic equations}
The first problem is related to the consideration of numbers of missing digits satisfying algebraic equations. Arithmetic properties of numbers with restricted digits have been studied for example in \cite{B15}, \cite{B13}, \cite{DMR11}, \cite{FM96}, \cite{MR10} as well as \cite{May2019}.  

Consider the equation
\[
x^3+y^3=1,
\]
and we want to ask whether or not there are irrational solutions $(x,y)$ such that both $x$ and $y$ do not have digit $1$ in their ternary expansions. More generally, we introduce the notion of  missing digits sets: Let $n\geq 1$ be an integer. Let $p\geq 3$ be an integer. Let $D\subset\{0,\dots,p-1\}^{n}.$ Consider the  set
\[
K_{p,D}=\cl{\{x\in\mathbb{R}^n: [p\{p^kx\}]\in  D,k\geq 0\}},
\] 
where $\{x\},[x]$ are the component wise fractional part, integer part respectively of $x\in\mathbb{R}^n$ and $\cl A$ is the closure of $A\subset\mathbb{R}^n$ under the standard topology.   More precisely, $\{x\}$ is the unique $y$ point in $[0,1)^n$ with $y-x\in\mathbb{Z}^n$, and $[x]=x-\{x\}\in\mathbb{Z}^n.$ 

For convenience, we use $\hat{D}$ for the complement $\{0,\dots,p-1\}^n\setminus D.$ For example, if $n=2,$ then $K_{3,\hat{\{(1,1)\}}}\cap [0,1]^2$ is the Sierpinski carpet with base $3$. Later on, we will call such a set $K_{p,D}$ to be a missing digits set. 

\begin{ques}\label{ques0}
	Let $n\geq 1$ be an integer. Let $M$ be a non-degenerate analytic manifold. Let $p>2$ be an integer and $D\subset \{0,\dots,p-1\}^n$ be a choice of digits. Determine whether or not $M\cap K_{p,D}$ is infinite.
\end{ques}

To guide our intuitions, we formulate the following perhaps ambitious conjecture.
\begin{conj}\label{conj1}
		Let $n\geq 1$ be an integer. Let $M$ be a strongly non-degenerate analytic manifold. Let $p>2$ be an integer and $D\subsetneq \{0,\dots,p-1\}^n$ be a choice of at least two digits. For $\delta>0,$ let $M^\delta$ be the $\delta$-neighbourhood of $M$ in $\mathbb{R}^n.$ Then $(M\cap K_{p,D})^\delta$ can be covered by
		\[
		\ll \left(\frac{1}{\delta}\right)^{\max\{0,\Haus K_{p,D}-(n-\Haus M)\}}
		\] 
		many $\delta$-balls.
\end{conj}
\begin{rem}
Strongly non-degeneracy is not a common notion and we will not use it anywhere in this paper. If $M$ is strongly non-degenerate, it means that $M$ is non-degenerate and for each affine subspace $L\subset\mathbb{R}^n,$
\[
\dim(M\cap L)\leq \max\{0,\dim M-(n-\dim L)\}.
\]
Intuitively, this condition says that $M$ is 'sufficiently twisted'. This condition is to avoid the cases when $K_{p,D}$ is contained in an affine subspace $L$ and $M\cap L$ has a larger dimension than expected so that $M\cap K_{p,D}$ can be larger than what is stated in the conjecture. Simple examples of strongly non-degenerate manifolds include unit spheres, Veronese curves, etc.
\end{rem}
\begin{rem}
	If the exponent of $(1/\delta)$ is equal to zero, this conjecture says that $M\cap K_{p,D}$ is a finite set. We can push a bit further. If in addition, $M$ is also an algebraic variety over $\mathbb{Q}$ and $\Haus K_{p,D}-(n-\Haus M)< 0$, then we expect that $M\cap K_{p,D}$ consists of only rational points.  For example, when $n=1$, this falls in the range of the aforementioned conjecture of \'{E}. Borel which says that all algebraic irrational numbers cannot miss digits in any base.
\end{rem}
\begin{rem}
We can formulate a slightly weaker conjecture with the covering number being
\[
\ll(1/\delta)^{\max\{0,\Haus K_{p,D}-(n-\Haus M)\}+\epsilon}
\]
for each $\epsilon>0$.
This weaker conjecture is also open in general. 
\end{rem}
One of the results we will prove is that Conjecture \ref{conj1} holds when $K_{p,D}$ is large enough. In that case, we also provide a natural lower counting estimate for missing digits points in $M$, see Theorem \ref{Main2}. It is of interest to find the 'smallest possible' missing digits set for the above conjecture to hold. As long as $M$ is fixed, we are able to provide a positive number $\sigma(M)>0$ and examples of $K_{p,D}$ with $\Haus K_{p,D}$ being larger than and arbitrarily close to $n-\sigma(M)$. Thus the missing digits sets need not be too large in the sense of Hausdorff dimension. In Section \ref{sec: Fail}, we demonstrate a particularly subtle difficulty in Conjecture \ref{conj1} for small missing digits sets.
\subsection{Intersecting manifolds with fractals}\label{sec: intersection}
We now discuss the number theoretic problem from a slightly different point of view. Let $n\geq 2$ and $M\subset \mathbb{R}^n$ be a manifold. Let $F\subset\mathbb{R}^n$ be a fractal set, e.g. a Sierpinski sponge. We are interested in considering the intersection $M\cap F.$ In view of the classical Marstrand slicing theorem (see \cite[Theorem 6.9 and Section 7]{Ma2}), we see that for a 'generic' translation vector $a\in\mathbb{R}^n$,
\begin{align*}\label{eqn: Dimension Reduction}
\Haus ((M+a)\cap F)\leq \max\{0,\Haus F-(n-\Haus M)\}.\tag{Dimension Reduction}
\end{align*}
Of course, it is possible to quantify the word 'generic' in a precise way. Loosely speaking, there is a small (in terms of the Lebesgue measure or Hausdorff dimension) exceptional set $E$ such that the above holds for all $a\notin E.$ In this direction, a great amount of efforts have been made to discover the occasions in which one can reduce the exceptional set from a small set to the empty set. For $M$ being affine subspaces, see \cite{Fu2}, \cite{HS12}, \cite{Sh}, \cite{Wu}.

Intuitively speaking, the only chance for the above (\ref{eqn: Dimension Reduction}) not to hold would be that $M$ and $F$ share similar structures in small scales. For example, if $M$ is a line parallel to one of the coordinate axis, then it is possible that \[
\Haus (M\cap F)>\max\{0,\Haus F-(n-\Haus M)\}.
\] 
This phenomenon happens in $\mathbb{R}^2$ already: just consider $F$ being the twofold Cartesian product  of the middle third Cantor set and $M$ being the $X$ coordinate line. In \cite{Sh}, Shmerkin showed that for $M$ being lines, those are essentially the only cases for (\ref{eqn: Dimension Reduction}) not to hold: (\ref{eqn: Dimension Reduction}) holds for all lines with irrational slopes.

Now, suppose that $M$ has some curved structures. Then intuitively, we can think that $M$ cannot share any structures with $F$ in small scales and (\ref{eqn: Dimension Reduction}) should hold without exceptions. 
\begin{comment}
This is indeed the case when $M$ is a hypersurface with non-vanishing curvatures. Indeed, with the help of  CP-chain method (introduced by Furstenberg in \cite{Fu2}), it is possible to prove the following result.
\begin{thm*}[CP-chain]
	Let $M\subset\mathbb{R}^n$ be a hypersurface with non-vanishing curvatures. Let $F$ be the Sierpinski sponge (based in 3,  say) with $\Haus F>n-1$. Then we have
	\[
	\Haus M\cap F\leq \max\{0,\Haus F-1\}.
	\]
\end{thm*} 

Hypersurfaces are submanifolds of codimension one. It is very natural to think about the question with submanifolds of any codimensions. 
\end{comment}
Towards this direction, we pose the following question.
\begin{ques}\label{ques1}
	Let $M\subset\mathbb{R}^n$ be a 'curved' submanifold. 
	Let $F$ be the Sierpinski sponge (based in 3,  say). Then we have
	\[
	\Haus M\cap F\leq \max\{0,\Haus F-n+\Haus M\}.
	\]
\end{ques}
Although both of the arguments in \cite{Wu} and \cite{Sh} do not work directly to gain information about intersections between  'curved' manifolds and fractals, it is perhaps possible to adapt the arguments to work out the case when the manifold is a hypersurface with non-vanishing curvatures.\footnote{We thank P. Shmerkin for explaining the details.} Here we do not fix the notion of a submanifold (possibly with codimension larger than one) being 'curved'. We leave the interpretation open. One of the results in this paper is to answer this question with a specific meaning to the word 'curved'. 
\subsection{Counting missing digits points near manifolds}
It turns out that the intersection problem above is closely related to a special lattice counting problem that we now introduce.

Recall the notion of missing digits sets $K_{p,D}.$ It is possible to introduce a natural Lebesgue measure $\lambda_{p,D}$ on $K_{p,D}$ which will we discuss later. For now, we just think $\lambda=\lambda_{p,D}$ to be a probability measure supported on $K=K_{p,D}.$

Let $M\subset\mathbb{R}^n$ be a submanifold. Let $\delta>0.$ Consider the $\delta$-neighbourhood $M^{\delta}$ of $M.$ We want to study the quantity
\[
\lambda(M^\delta)
\]
for $\delta\to 0$. Heuristically, if $\lambda$ and $M$ are somehow 'independent', we expect that
\begin{align*}\label{eqn: Independence}
\lambda(M^\delta)\asymp\delta^{n-\dim M}.\tag{Independence}
\end{align*}
Here $\dim M$ is the standard notion of the dimension of the submanifold $M$. In our situation, $\dim M=\Haus M.$ The value $n-\dim M$ is usually called the codimension of $M.$  The value $\delta^{n-\dim M}$ is roughly the Lebesgue measure of $M^\delta.$  Now, assuming (\ref{eqn: Independence}), we see that $M^\delta\cap K$ can be covered with
\[
\ll\delta^{n-\dim M}/\delta^{\Haus K}
\]
many balls with radius $\delta.$ For this, we need to assume that $\lambda$ is AD-regular with exponent $\Haus K.$ See Section \ref{sec: AD}. From here, we directly deduce that
\[
\Haus (M\cap K)\leq \Haus K-(n-\dim M)
\]
if $\Haus K>n-\dim M.$ Otherwise, we simply have
\begin{align*}\label{eqn: Finiteness}
\# M\cap K<\infty.\tag{Finiteness}
\end{align*}
The conclusion (\ref{eqn: Finiteness}) is rather strong, it says that if $K$ is a small enough missing digits set then $K\cap M$ has only finitely many points. Of course, we would have to assume the asymptotic bound $\lambda(M^\delta)\asymp \delta^{n-\dim M}$ which is not easy to be tested. A particular special case (of Question \ref{ques0}) in this direction can be formulated as follows.
\begin{ques*}
Consider the circle $x^2+y^2=1$ in $\mathbb{R}^2.$ Can we find infinitely many points $(x,y)$ on the circle with $x,y\in K_{5,\{0,4\}}$?
\end{ques*}
This is an open problem. Methods in \cite{Sh} or \cite{Wu} can be probably used to deduce that the points under consideration form a set with zero Hausdorff dimension but this is not enough to deduce the finiteness. More generally, we shall consider the following problem.
\begin{ques}\label{ques2}
	Find examples of $M$ and $K_{p,D}$ (and $\lambda_{p,D}$) for which the asymptotic estimate holds
	\[
	\lambda_{p,D}(M^\delta)\asymp \delta^{n-\dim M}.
	\]
\end{ques}

\section{results in this paper}
To state the results, we first introduce some terminologies. The reader can skip to Section \ref{results} and return here later when it is necessary.
\subsection{Fourier norm dimensions }\label{sec: Fourier Dimension}
We need some notions of Fourier norm dimensions. They are useful in e.g. \cite{Y21} where a problem of counting rational points near missing digits sets was considered. In an ongoing project, \cite{ACVY}, an algorithm for computing the Fourier $l^1$ dimensions of missing digits sets is developed together with many applications to metric Diophantine approximation. In this paper, we do not need to consider the precise values of the Fourier $l^1$-dimensions.  We only provide some bounds which are rather crude but enough for the discussions in this paper. See Theorem \ref{thm: l1 bound}.

Let $n\geq 1$ be an integer. Let $\mu$ be a compactly supported Borel probability measure on $\mathbb{R}^n.$ Consider the Fourier transform
\[
\hat{\mu}(\xi)=\int_{\mathbb{R}^n} e^{-2\pi i (x,\xi)}d\mu(x),
\]
where $(.,.)$ is the standard Euclidean bilinear form. 
\begin{defn}
Let $p>0.$ We define
	\[
	\dim_{l^p}\mu=\sup\left\{s>0: \sup_{\theta\in [0,1]^n}\sum_{|\xi|\leq R,\xi\in\mathbb{Z}^n} |\hat{\mu}(\xi+\theta)|^p \ll R^{n-s}\right\}.
	\] 
\end{defn}
With the help of the Cauchy-Schwarz inequality, it is possible to show that
\[
\frac{\dim_{l^2}\mu}{2}\leq \dim_{l^1}\mu.
\]
Moreover, we have for each AD-regular (see Section \ref{sec: AD}) measure $\mu$
\[
\dim_{l^2}\mu=\Haus \mu=\Haus \supp(\mu).
\]
Furthermore, let $n\geq 1$ be an integer. Let $\mu_1,\dots,\mu_n$ be a Borel probability measure on $\mathbb{R}.$ The $n$-fold Cartesian product $\mu'=\mu_1\times\dots\times \mu_n$ satisfies
\[
\dim_{l^1}\mu'\geq \dim_{l^1}\mu_1+\dots+\dim_{l^1}\mu_n.
\]
In fact, the equality holds when $\mu_1,\dots,\mu_n$ are missing digits measures but we do not need this fact.

We have seen above that $\dim_{l^2}\mu$ is closely related to $\Haus \mu.$ The reason that we also study $\dim_{l^1}\mu$ is that it gauges, in some sense, how 'close' is $\mu$ from being a continuous function. Observe that if the exponent in the definition of $\dim_{l^1}\mu$ can be chosen to be negative, then $\mu$ has an absolutely integrable Fourier transform. This says that $\mu$ can be chosen to be the distribution associated with a continuous density function. In this case, $\supp \mu$ can be seen as a topological manifold. 

For computations, it is often not very convenient to have two $\sup$'s. For this reason, we also introduce the following two additional definitions.
\begin{defn}[Integral]
	Let $p>0.$ We define
	\[
	\dim^{I}_{l^p}\mu=\sup\left\{s>0: \int_{|\xi|\leq R} |\hat{\mu}(\xi)|^pd\xi\ll R^{n-s}\right\}.
	\] 
\end{defn}
\begin{defn}[Sum]
Let $p>0.$ We define
	\[
	\dim^{S}_{l^p}\mu=\sup\left\{s>0: \sum_{|\xi|\leq R,\xi\in\mathbb{Z}^n} |\hat{\mu}(\xi)|^p \ll R^{n-s}\right\}.
	\] 
\end{defn}
In most of the situations, $\dim^I_{l^p}, \dim^S_{l^p}$ can provide useful information individually. Notice that in general,
\[
\dim_{l^p}\mu\leq \min\{\dim^I_{l^p}\mu,\dim^S_{l^p}\mu\}.
\]
This is because
\begin{align*}
&\int_{|\xi|\leq R} |\hat{\mu}(\xi)|^pd\xi\leq \sum_{\xi\in\mathbb{Z}^n, |\xi|\leq 1.5R}\int_{\theta\in [0,1]^n} |\hat{\mu}(\xi+\theta)|^pd\theta\\
&=\int_{\theta\in [0,1]^n} \sum_{\xi\in\mathbb{Z}^n, |\xi|\leq 1.5R}|\hat{\mu}(\xi+\theta)|^pd\xi\\
&\leq \sup_{\theta\in [0,1]^n}\sum_{|\xi|\leq 2R,\xi\in\mathbb{Z}^n} |\hat{\mu}(\xi+\theta)|^p
\end{align*}
and
\[
 \sum_{|\xi|\leq R,\xi\in\mathbb{Z}^n} |\hat{\mu}(\xi)|^p\leq \sup_{\theta\in [0,1]^n}\sum_{|\xi|\leq R,\xi\in\mathbb{Z}^n} |\hat{\mu}(\xi+\theta)|^p.
\]
We also suspect that for missing digits measures, the three notions $$\dim_{l^p},\dim^I_{l^p},\dim^S_{l^p}$$ are identical.
\subsection{Manifolds of finite type}
Let $n\geq 2$ be an integer. A smooth submanifold $M\subset \mathbb{R}^n$ is of finite type if for each $x\in M,$ the manifold $M$ only has finite order contacts with affine hyperplanes at $x$. For more details, see \cite[Chapter VIII, Section 3]{Stein}. We do not repeat the definition here and only illustrate some of the examples.

(1). Let $M$ be an analytic submanifold. Suppose that $M$ is not contained in any affine hyperplane, then $M$ is of finite type.

(2). As a particular example, consider the Veronese curve $(t,t^2,t^3,\dots,t^n),t\in\mathbb{R}.$ This curve is analytic and it is not contained in any affine hyperplane. Therefore it is of finite type.

(3) If $M$ is a smooth hypersurface with non-vanishing Gaussian curvature, then $M$ is of finite type.
\begin{comment}
(3). More generally, irreducible algebraic varieties are almost manifolds of finite type. There are some complications as real varieties are not always manifolds. Sometimes, there are complicated singularities. However, small enough neighbourhoods around non-singular points on varieties are analytic manifolds. Thus locally one can apply the criteria from (1).
\end{comment}
\subsection{Missing digits sets and measures}
We recall the notion of missing digits sets. Let $n\geq 1$ be an integer. Let $p\geq 3$ be an integer. Let $D\subset\{0,\dots,p-1\}^{n}.$ Consider the  set
\[
K_{p,D}=\cl\{x\in\mathbb{R}^n: [p\{p^kx\}]\in  D,k\geq 0\},
\] 
where $\{x\},[x]$ are the component wise fractional part, integer part respectively of $x\in\mathbb{R}^n.$
Let $p_1,\dots,p_{\#D}$ be a probability vector, i.e. they are non-negative and sum to one. We can then assign each element in $D$ a probability weight. To be specific, one can first introduce an ordering on $D$ and assign the probabilities accordingly.  We can now construct the random sum
\[
S=\sum_{i\geq 1} p^{-i} \mathbf{d}_i
\]
where $\mathbf{d}_i\in D,i\geq 1$ are randomly and independently chosen from the set $D$ with the assigned probabilities.

If $p_1=\dots=p_{\#D}=1/\#D,$ the distribution of $S$ is a Borel probability measure supported on $[0,1]^n.$ We call this measure to be $\lambda_{p,D}.$ It is a Borel probability measure supported on $K_{p,D}\cap [0,1]^n.$ Moreover, it is AD-regular with exponent $\Haus K_{p,D}.$ We also write
\[
\dim_{l^1}K_{p,D}=\dim_{l^1}\lambda_{p,D}.
\]

We provide some examples of missing digits sets. Recall that $\hat{D}=\{0,\dots,p-1\}^n\setminus D.$

(1) If $n=1,$ then $K_{3,\{0,2\}}\cap [0,1]$ is the middle third Cantor set and $\lambda_{3,\{0,2\}}$ is the natural middle third Cantor-Lebesgue measure.

(2) If $n=2,$ then $K_{3,\hat{\{(1,1)\}}}\cap [0,1]^2$ is the Sierpinski Carpet with base $3.$ In general, if $n\geq 3,$ $K_{3,\hat{\{(1,\dots,1)\}}}\cap [0,1]^n$ are Sierpinski sponges.

\begin{rem}
It is often interesting to consider Cartesian products of missing-digit sets (measures) which are not necessarily missing-digit sets (measures) themselves. In fact, they are not self-similar in general. Although we do only consider missing-digit sets (measures) in this paper, many results can be extended to deal with Cartesian products of missing-digit sets (measures). We will discuss this in Section \ref{sec: application}. Before that, one can ignore this technical consideration.
\end{rem}
\subsection{Results}\label{results}
Towards Conjecture \ref{conj1}, we prove the following theorem. Later on, we will provide examples that fall into the range of this theorem. The conditions are related to the $l^1$-dimensions of missing digits sets. This is not a common notion of dimension. Loosely speaking, missing digits sets with large bases have almost equal Hausdorff and $l^1$-dimensions. So it is helpful to think $\dim_{l^1}$ below just as $\Haus.$\footnote{In Theorem \ref{thm: l1 bound}(2) we will see a precise version of this 'loose' statement.}
\begin{thm}\label{Main}
	Let $n\geq 2$ be an integer. Let $M\subset\mathbb{R}^n$ be a manifold of finite type. Then there is a number $\sigma=\sigma(M)>0$ such that for each compactly supported Borel probability measure $\lambda$ with $\dim_{l^1}\lambda>n-\sigma,$
	\[
	\lambda(M^\delta)\ll \delta^{n-\dim M}.
	\]
\end{thm}
This proves Conjecture \ref{conj1} for large missing digits sets.  The number $\sigma(M)$ can be explicitly determined once we know $M.$ It is related to the Fourier decay properties of smooth surface measures carried by $M.$ This number is always $\leq \dim M/2.$ In the case when $M$ is a hypersurface with non-vanishing Gaussian curvature, it can be chosen to be $\dim M/2=(n-1)/2.$  The condition we have in this theorem is in some sense sharp. We postpone the discussion to Theorem \ref{thm: fail}.

We will provide some crude bounds for $\dim_{l^1}\lambda_{p,D}$, see Theorem \ref{thm: l1 bound}. In particular, the following result can be deduced.
\begin{cor}[Missing digits points near manifolds]\label{cor: near manifold}
	Let $n\geq 2$ be an integer. Let $M\subset\mathbb{R}^n$ be a manifold of finite type. Let $k\geq 1$ be an integer. There is a number $p_0(M,k)$ so that for each integer $p>p_0(M,k)$,  Theorem \ref{Main} holds for $\lambda_{p,D}$ where the digit set $D$ satisfies $\#D\geq p^n-k.$
\end{cor}
The number $p_0(M,k)$ can be explicitly determined once we know $\sigma(M)$ and $k.$ Theorem \ref{Main} leads to an intersection result.
\begin{cor}[Intersections]\label{cor: intersection}
	Let $n\geq 2$ be an integer. Let $M\subset\mathbb{R}^n$ be a manifold of finite type. Then there is a number $\sigma=\sigma(M)>0$ such that for each compactly supported Borel probability measure $\mu$ which is AD-regular and $\dim_{l^1}\mu>n-\sigma,$
	\[
	\ubox M\cap supp(\mu)\leq \max\{0,\Haus \supp(\mu)-(n-\dim M)\}.
	\]
\end{cor}
On the other hand, we are also interested in having a lower estimate for $M\cap K_{p,D}.$ However, this problem is not always sensible to ask because $K_{p,D}$ has holes and it can happen that $M$ is contained in one of the holes. This phenomenon is certainly not capturing the real nature between $M$ and $K.$ In order to unveil the underlying truth. We propose two compromises by modifying $K$ or $M$ in 'acceptable' manners that will be clarified soon.

First, we can render this problem by relaxing our requirement for missing digits sets. Let $l\geq 0$ be an integer. We define
\[
K_{p,D,l}=\cl\{x\in\mathbb{R}^n: [p\{p^kx\}]\in D, k\geq l\}.
\]
So $K_{p,D,l}$ fills some of the holes of $K_{p,D}.$ The purpose of doing this is to get rid of the special cases when the manifold just happens to be contained in the holes of the missing digits set. We do not relax the situation too much by requiring that $l$ has to be a fixed integer. For example, the digital representations of $x\in K_{p,D,l}\cap [0,1]^n$ are of no restriction only for the first finitely many digits. This is the 'acceptable' modification for $K$ we mentioned earlier.
\begin{thm}\label{Main2}[Compromise 1]
	Let $M$ be a manifold of finite type. Let $K=K_{p,D}$ be a missing digits set. Suppose that $\dim_{l^1} K_{p,D}+\sigma(M)-n>0$. Then there is an $l\geq 0$ such that the set
	\[
	M\cap K_{p,D,l}
	\]
	is infinite. In fact, there is an $l\geq 0$ such that for all small enough $\delta>0,$ to cover $M^\delta\cap K_{p,D,l}^\delta,$ at least $\gg (1/\delta)^{\Haus K+\dim M-n}$ many $\delta$-balls are needed.  Moreover, we have
	\[
	\lbox M\cap K_{p,D,l}\geq \Haus K-(n-\sigma(M)).
	\]
	Finally,  $\bigcup_{l\geq 0} K_{p,D,l}\cap M$ is dense in $M.$
\end{thm}
\begin{rem}
	It is tempting to believe that $$\boxd K_{p,D,l}\cap M=\Haus K-(n-\dim M).$$ So far, we can conclude that
	\[
	\Haus K-n+\sigma(M)\leq\lbox M\cap K_{p,D,l}\leq \ubox M\cap K_{p,D,l}\leq \Haus K-(n-\dim M).
	\]
	Since $\sigma(M)\leq \dim M/2,$ the above chain of inequalities will never be close to being optimal. Some other ideas need to be used in order to fill in this gap.
	
	The last statement seems trivial. The set $\cup_{l\geq 0} K_{p,D,l}$ is dense in $\mathbb{R}^n.$ However, without the previous statement, it is even unclear whether or not $\cup_{l\geq 0} M\cap K_{p,D,l}$ is empty.
	\end{rem}
	
	\begin{comment}

\begin{thm}\label{lower bound}
Let $M$ be a non-degenerate analytic hypersurface. Let $K$ be a product of missing digits sets on $\mathbb{R}$ with coprime bases and $\dim_{l^1}K>(n+1)/2.$ Suppose that
\[
K=K_{p_1,D_1}\times\dots\times K_{p_n,D_n}
\]
where $(p_1,\dots,p_n)=1$ and $D_1,\dots,D_n$ are proper digit sets. Then there are integers $l_1,\dots,l_n\geq 0$ such that for
\[
K'=K_{p_1,D_1,l_1}\times\dots\times K_{p_n,D_n,l_n}
\]
we have
\[
\boxd M\cap K'=\Haus K-(n-(n-1))=\Haus K-1.
\]
\end{thm}
\end{comment}
Another possible way to compromise the fact that $K=K_{p,D}$ has holes is to consider a family of transformed manifolds of $M.$ In this paper, we deal with the group $E_n=(0,\infty)\times \mathbb{R}^n\times\mathbb{O}(n).$ For $(t,v,g)\in E_n,$ it acts on $\mathbb{R}^n$ by first applying  the rotation $g$, then applying the scaling $t$, then applying the translation $v$:
\[
x\in\mathbb{R}^n\to T_{t,v,g}(x)=tg(x)+v.
\]
Observe that $E_n$ also acts on the space of non-degenerate manifolds and it will not distort Fourier transforms too much. Thus we expect that
\begin{align*}\label{Dimension Conservation}
\boxd T_{t,v,g}(M)\cap K=\Haus K-(n-\dim M)\tag{Dimension Conservation}
\end{align*}
and in particular $T_{t,v,g}(M)\cap K\neq\emptyset$
for many choices of $(t,v,g)\in E_n.$ In view of Mattila's intersection theorem (\cite[Section 7]{Ma2}), we already know that the above holds for 'many' choices of $(t,v,g)$ in a metric sense. We now upgrade the metric result at a cost of dropping the (\ref{Dimension Conservation}).
\begin{thm}\label{Main3}[Compromise 2]
    Let $n\geq 2$ be an integer. Let $M$ be a manifold of finite type. Suppose further that $\lambda$ is a missing digits measure supported on $K.$ The function
    \[
    f:(t,v,g)\to \limsup_{\delta\to 0} \frac{\lambda((T_{t,v,g}(M))^\delta)}{\delta^{n-\dim M}}
    \]
    is well-defined and non-negative valued. If $M$ is compact, then there are a constant $c>1$ and a continuous and non-vanishing function $h: E_n\to [0,\infty)$ such that
    \[
    c^{-1}h\leq f\leq ch.
    \]
\end{thm}
\begin{rem}\label{Remark: unproved}
(a). The compactness of $M$ in the last statement is not strictly necessary. In fact, if $M$ is not compact, one can replace $M$ with a compactly supported and smooth surface measure on $M$.

(b). This result is likely to hold for other classes of possibly nonlinear transformations of manifolds other than the scalings, translations and rotations, for example, the evolution of $M$ under a smooth vector field with some regularity conditions. As $E_n$ is already enough for most of our arithmetic applications we do not pursue this degree of generality.\footnote{This paragraph is not mathematically precise. We only want to illustrate an idea rather than a rigorous definition. One possible way of seeing the idea in this paragraph is to consider a light and soft tissue floating in the air.}

(c). We are curious about whether or not the function $f$ can be defined via $\mathcal{H}^s$ for $s=\Haus K-(n-\dim M).$ More precisely, let
\[
h:(t,v,g)\to \mathcal{H}^s(T_{t,v,g}(M)\cap K).
\]
Is $h$ continuous? From Theorem \ref{Main}, it is not hard to show that $h$ takes real values (not $\infty$). In this way, one can the understand the Hausdorff dimension of $T_{t,v,g}(M)\cap K.$ 

(d). Notice that whenever $f(t,v,g)>0,$ it follows that
\[
T_{t,v,g}(M)\cap K\neq\emptyset.
\]
However, this is not enough to conclude that $\boxd$ or $\Haus T_{t,v,g}(M)\cap K\geq s.$ In this case, we have the weaker result
\[
\lbox T_{t,v,g}(M)\cap K\geq \Haus K-(n-\sigma(M))
\]
as in Theorem \ref{Main2}.
\end{rem}
The above results also hold when $M$ is replaced with other sets supporting measures with polynomial Fourier decay. We will discuss this part in Section \ref{sec: fractal measure}.

Finally, we emphasize that the condition $\dim_{l^1}\lambda>n-\sigma(M)$ cannot be dropped although we believe that it might not be the optimal condition to look at. See Section \ref{sec: further}.
\begin{thm}[Sharpness]\label{thm: fail}
Let $M\subset\mathbb{R}^2$ be the curve
\[
M: |x|^k+|y-1|^k=1
\]
where $k\geq 2$ is an integer. Let $p>3$ be an integer and let $D_1=\{0,\dots,l_1\}$, $D_2=\{0,\dots,l_2\}$ for some integers $0<l_1,l_2<p-1.$ Let $K=K_{p,D_1}\times K_{p,D_2}\cap [0,1]^2.$ Then for $\delta>0,$
$
M^\delta\cap K
$
cannot be covered with
\[
\ll \left(\frac{1}{\delta}\right)^{(k-1)s/k}
\]
many $\delta$-balls where $s=\max\{\Haus K_{p,D_1}, \Haus K_{p,D_2}\}.$ 
\end{thm}
\begin{rem}\label{rem: fail}
In particular, if $k=2,$ then $M$ is a circle. In this case, the threshold for $\dim_{l^1}$ in Theorem \ref{Main} is $n-\sigma(M)=2-1/2=3/2.$ It is possible to choose $D_1,D_2$ such that $\Haus K_{p,D_1}$ is close to one and $\Haus K_{p,D_2}$ is close to $1/2.$ Then $\Haus K_{p,D_1}\times K_{p,D_2}$ is close to $3/2.$ Moreover, if $p$ is large enough, then $\dim_{l^1}K_{p,D_1}\times K_{p,D_2}$ is also close to $3/2$. See Theorem \ref{thm: l1 bound}. In this case we see that $M^\delta\cap K$ cannot be covered with
\[
\ll \left(\frac{1}{\delta}\right)^{s/2}
\]
many $\delta$-balls. Notice that $s/2$ can be made arbitrarily close to $1/2$. On the other hand, if the conclusion of Theorem \ref{Main} would hold, then $M^\delta\cap K$ can be covered by
\[
\ll \left(\frac{1}{\delta}\right)^{\Haus K-1}
\]
many $\delta$-balls. Therefore if $\Haus K-1<s/2,$ then the conclusion of Theorem \ref{Main} cannot hold. Since $s$ can be chosen to be close to one, we are able to find examples of missing digits measures $\lambda$ with $\dim_{l^1}\lambda$ smaller but arbitrarily close $3/2$ such that the conclusion of Theorem \ref{Main} does not hold for $\lambda$ and $M.$

In general, the above discussion works for $k\geq 3$ as well. In this case, $M$ is a 'flatter' curve than the circle. We have $n-\sigma(M)=2-(1/k).$ This is because $M$ has $k$th order contact with the $X$-axis. See Section \ref{sec: Fail} for more details. On the other hand, as above we see that if $\Haus K-1<(k-1)s/k,$ then Theorem \ref{Main} cannot hold for the missing digit measure on $K.$ Thus we are able to find examples of missing digits measures $\lambda$ with $\dim_{l^1}\lambda$ smaller but close to $2-(1/k)$ such that the conclusion of Theorem \ref{Main} does not hold for $\lambda$ and $M.$ 
\end{rem}
\begin{rem}
Theorem \ref{thm: fail} does not disprove Conjecture \ref{conj: Palis}. However, it does show that it is perhaps a challenging task to reduce the $3/4$ threshold. Theorem \ref{thm: fail} does not disprove Conjecture \ref{conj1} either, because in the statement of Conjecture \ref{conj1}, what we are interested in is to cover the set $(M\cap K)^\delta,$ which is in general smaller than $M^\delta\cap K^\delta.$
\end{rem}
\subsection{Fourier norm dimensions of missing digits measures}
Missing digits sets and measures provide us with a large class of examples of fractal sets and measures for the results in the previous section to be applied. In this section, we list a few general results regarding the $l^1$-dimensions of missing digits sets whose proofs are provided to make this paper self-contained\footnote{In \cite{ACVY}, a much more precise method is developed.}. 

\begin{thm}\label{thm: l1 bound}
	Let $n\geq 1$ be an integer. The following results hold.
	\begin{itemize}
		\item{1} Let $t\geq 1$ be an integer. We have
		\[
		\liminf_{p\to\infty,\#D\geq p^n-t} \dim_{l^1} \lambda_{p,D}=n.
		\]
		In particular, for each number $\epsilon>0,$ as long as $p$ is large enough, $\dim_{l^1}\lambda_{p,D}>n-\epsilon$ holds for each $D$ with $\#D=p^n-1.$
		\item{2} For each integer $p\geq 4,$ let $D\subset\{0,\dots,p-1\}^n$ be a 'rectangle', i.e. a set of form $[a_1,b_1]\times [a_2,b_2]\dots [a_n,b_n]\cap \{0,\dots,p-1\}^n.$ Then we have \footnote{The base of $\log$ in this paper is $e$.}
		\begin{align*}\label{T}
		\dim_{l^1}\lambda_{p,D}\geq \Haus \lambda_{p,D}-\frac{n\log\log p^2}{\log p}.\tag{T}
		\end{align*}
		The significance of this result is to show that missing digits sets can have almost equal Hausdorff and $l^1$-dimensions.
	\end{itemize}
\end{thm}
\begin{rem}
The digit set structure is rather special for part 2. We emphasize that it can be a little bit more complicated. More precisely, part 2 of this theorem also works for rectangles with gaps larger than $1.$ For example, instead of being a product set of sets of consecutive integers, it can be a product set of sets of consecutive even numbers. Also, one can consider unions of rectangles as long as the number of rectangles in the union is not too large compared with $p,$ e.g. $O(p^{\epsilon})$. We are interested in finding the most general condition for the digit sets so that a lower bound like (\ref{T}) holds, i.e. as $p\to\infty,$
\[
|\dim_{l^1}\lambda_{p,D}-\Haus \lambda_{p,D}|=o(1).
\]
\end{rem}
Towards the question in this remark, we pose the following possibly too optimistic conjecture.
\begin{conj}\label{conj: l1 bound}
	Let $p>2$ be a prime number. Let $\epsilon\in (0,1).$  Then we have
	\[
	|\dim_{l^1}\lambda_{p,D}-\Haus \lambda_{p,D}|=o(1),
	\]
	where the $o(1)$ term is uniform across all digit set $D$ with
	\[
	\# D\geq p^{\epsilon}.
	\]
\end{conj}
The condition that $p$ is a prime number cannot be completely dropped. In fact, if $p=3^n$ where $n$ is an integer. Then it is possible to choose $D$ such that $K_{p,D}$ is the middle third Cantor set. The statement of this conjecture is obviously not true in this case.

\section{Other related works}
The problem in this paper can be generally described as "understanding the distribution of a special class of points near another set." Consider the following three classes of sets:

1. Manifolds

2. Rational numbers with bounded heights (denominators).

3. Missing digits sets.

We studied the distribution of 3 near 1 in this paper. On the other hand, it has been of great interest in considering the distribution of 2 near 1 (see \cite{B12}, \cite{BVVZ}, \cite{BDV}, \cite{KM98}, \cite{H20}, \cite{VV06}) as well as 2 near 3 (see \cite{ACY}, \cite{ACVY}, \cite{BFR11}, \cite{BD16}, \cite{EFS11}, \cite{KL20}, \cite{KLW}, \cite{LSV ref}, \cite{PV2005}, \cite{SW19}, \cite{W01}, \cite{Y21}).\footnote{Those are just a small fraction of all the works on those topics.} 

The Fourier analysis method in this paper can be found in many of the above references.  Other than this Fourier analysis method, another major tool for considering the above counting problems is by using the theory of dynamics on homogeneous spaces, e.g. \cite{KM98} (2 near 1), \cite{KL20} (2 near 3). A dynamical approach for the problems considered in this paper (3 near 1) is likely to be found. 
\section{Preliminaries}
Before we prove Theorem \ref{Main}, we need to introduce some more definitions as well as some results which will be used without proofs. See \cite{Fa}, \cite{Ma1} for more details on the notions of dimensions and regularities of measures.
\subsection{Hausdorff dimension, box dimensions}
Let $n\geq 1$ be an integer. Let $F\subset\mathbb{R}^n$ be a Borel set. Let $g: [0,1)\to [0,\infty)$ be a continuous function such that $g(0)=0$. Then for all $\delta>0$ we define the  quantity
\[
\mathcal{H}^g_\delta(F)=\inf\left\{\sum_{i=1}^{\infty}g(\mathrm{diam} (U_i)): \bigcup_i U_i\supset F, \mathrm{diam}(U_i)<\delta\right\}.
\]
The $g$-Hausdorff measure of $F$ is
\[
\mathcal{H}^g(F)=\lim_{\delta\to 0} \mathcal{H}^g_{\delta}(F).
\]
When $g(x)=x^s$ then $\mathcal{H}^g=\mathcal{H}^s$ is the $s$-Hausdorff measure and Hausdorff dimension of $F$ is
\[
\Haus F=\inf\{s\geq 0:\mathcal{H}^s(F)=0\}=\sup\{s\geq 0: \mathcal{H}^s(F)=\infty          \}.
\]

Next, let $F\subset\mathbb{R}^n$ be a Borel set. Let $\delta>0$ and $N_{\delta}(F)$ be the minimum amount of $\delta$-balls needed to cover $F.$ Then the upper box dimension of $F$ is
\[
\ubox F=\limsup_{\delta\to\infty} \frac{-\log N_\delta(F)}{\log \delta}.
\]
The lower box dimension is
\[
\lbox F=\liminf_{\delta\to\infty} \frac{-\log N_\delta(F)}{\log \delta}.
\]
If the upper and lower box dimensions of $F$ are equal, we call the comment value to be the box dimension of $F$
\[
\boxd F=\lbox F=\ubox F.
\]
If $F$ is compact, then we have the general result
\[
\Haus F\leq \lbox F\leq \ubox F.
\]
\subsection{AD-regularity}\label{sec: AD}
Let $n\geq 1$ be an integer. Let $\mu$ be a Borel measure. Let $s>0$ be a number. We say that $\mu$ is $s $-regular, or AD-regular with exponent $s$ if there is a constant $C>1$ such that for all $x\in supp(\mu)$ and all small enough $r>0$
\[
C^{-1} r^s\leq \mu(B_r(x))\leq C r^s,
\] 
where $B_r(x)$ is the Euclidean ball of radius $r$ and centre $x$. For an AD-regular measure $\mu$, the exponent can be seen as
\[
s=\Haus supp(\mu).
\]
For this reason, we simply define for AD-regular measure $\mu,$
\[
\Haus \mu=\Haus supp(\mu).
\]
Missing digits measure $\lambda_{p,D}$ in $\mathbb{R}^n$ are AD-regular measures with exponent
\[
s=\Haus \lambda_{p,D}=\Haus K_{p,D}=\frac{\log \#D}{\log p^n}.
\]
\subsection{Fourier transform of surface measures on manifolds of finite type}\label{sec: fourier manifold}
Discussions in this section follow \cite[Chapter VIII, Section 3]{Stein}.

Let $n\geq 2$ be an integer and $M\subset\mathbb{R}^n$ be a manifold of finite type. Let $\mu$ be the surface measure, i.e. the natural Lebesgue measure carried by the manifold. If $M$ is compact, we can normalize $\mu$ to be a probability measure. Otherwise, we shall truncate the measure to vanish at infinity by choosing a smooth, compactly supported function $\phi$ and denote
\[
\mu'=c\phi\mu
\]
where $c$ is the normalization factor so that $\mu'$ is a probability measure.

We will always perform the truncating procedure. This will cause no loss of generality for the problems considered in this paper. 

A standard result of smooth surface measures carried by manifolds of finite type is that the Fourier transforms decays polynomially,
\[
\hat{\phi\mu}(\xi)=O(|\xi|^{-\sigma}).
\] 
Where $\phi$ is a smooth, compactly supported function on $\mathbb{R}^n.$ Of course, if $M$ is compact, one can also let $\phi=1.$

Here lower bounds for $\sigma$ can be effectively determined. It is related to the type of $M.$ Roughly speaking, $\sigma=1/k$ where $k$ is the smallest integer such that $M$ does not have $k$-th order contact with affine hyperplanes.  For the Veronese curve, one can choose $\sigma=1/n.$ If $M$ is a hypersurface with non-vanishing curvatures, then the Fourier transform has a much better decay: one can choose $\sigma=(n-1)/2=2^{-1}\dim M.$

The choice of $\sigma(M)$ for non-degenerate manifolds is a challenging topic in harmonic analysis and differential geometry. Apart from some general results in Stein's book \cite{Stein}, the article \cite{CI} contains some more technical results that are useful to provide estimates of $\sigma(M)$ for some specially defined $M$, e.g. joint graphs of analytic maps. 
\subsection{Fourier transform of missing digits measures}
The discussion in this subsection works for all self-similar measures with a uniform contraction ratio. We nonetheless focus only on missing digits measures. 

Let $n\geq 1$ be an integer. Let $p>2$ be an integer and $D\subset\{0,\dots,p-1\}^n.$ Consider the missing digit measure $\lambda_{p,D}.$ In this case, $\hat{\lambda}_{p,D}(\xi)$ can be computed with the help of the formula
\[
\hat{\lambda}_{p,D}(\xi)=\prod_{j\geq 0} \frac{1}{\#D}\sum_{d\in D} e^{-2\pi i (d,\xi)/p^j}.
\]
\subsection{Asymptotic notations}
We use both the Vinogradov ($\ll,\gg,\asymp$) as well as Bachmann-Landau ($O(),o()$) notations:

Let $f(\delta), g(\delta)$ be two real valued quantities depending on $\delta>0.$ Then
\begin{itemize}
    \item $f\ll g$ or $f=O(g)$ if $|f(\delta)|\leq C|g(\delta)|$ for a constant $C>0$ and all $\delta>0.$
    
    \item $f= o(g)$ if for each $\epsilon>0,$ there is a $\delta_0>0$ such that for all $\delta<\delta_0,$
    \[
    |f(\delta)|\leq \epsilon |g(\delta)|.
    \]
    
    \item $f\asymp g$ if $f\ll g$ and $g\ll f$.
\end{itemize}

\section{Proof of the results}
\subsection{Thickening the surface measure}
Let $n\geq 2$ be an integer. Let $M\subset\mathbb{R}^n$ be a smooth submanifold of finite type. Let $\mu$ be a probability measure supported on $M.$ We assume that $\mu$ is compactly supported and smooth. Those conditions are not necessary for the following arguments, we only assume them for convenience.

For later considerations, we need to thicken the measure $\mu.$ Let $\delta>0.$ If $\delta$ is small enough, the neighbourhood $M^\delta$ containing points in $\mathbb{R}^n$ that are $\delta$ close to $M$ provides us with a good approximation of $M.$  Let $\phi$ be a compactly supported and smooth function which equals to one on $B_{1/2}(0)$ and vanish outside of $B_{1}(0).$ We also arrange that $\phi$ and $\hat{\phi}$ are  spherically symmetric positive valued functions. See \cite{Y20b}.  Furthermore, we shall assume that $\int \phi=1.$ Let $\phi_\delta(x)=\delta^{-n}\phi(\delta x).$

We see that the convoluted measure $\mu_\delta=\phi_\delta*\mu$ is a probability measure and it is also a smooth function supported on $M^\delta.$ It can be also checked that on $M^{\delta/2},$ $\mu_\delta\asymp \delta^{-(n-\dim M)}$ and in general
\[
\mu_\delta(x)\ll \delta^{-(n-\dim M)}
\]
uniformly across $x\in\mathbb{R}^n.$

Observe that
\[
\hat{\phi_\delta*\mu}(\xi)=\hat{\phi}_\delta(\xi)\hat{\mu}(\xi).
\]
The function $\hat{\phi}$ decays fast at infinity. In fact, for each number $N>0,$ we have
\[
|\hat{\phi}(\xi)|= O(|\xi|^{-N}).
\]
Intuitively (after scaling by $\delta$), this tells that $\hat{\phi}_\delta$ is a smooth function which is essentially supported on $B_{1/\delta}(0).$ Within the ball $B_\delta(0)$ we have
\[
{\phi}_\delta\asymp \delta^{-n},
\]
that is, the value of $\phi_\delta/\delta^{-n}$ on $B_{\delta}(0)$ is bounded and strictly positive (in a manner than does not depend on $\delta$). Similarly, on $B_{1/\delta}(0),$ we have
\[
\hat{\phi}_\delta\asymp 1.
\]
\subsection{Littlewood-Paley decomposition, proof of Theorem \ref{Main}}\label{sec: L2}
Let $\lambda$ be a Borel probability measure, by Plancherel's theorem, we have
\[
\int \mu_\delta(x)d\lambda(x)=\int \hat{\mu}_\delta(\xi)\overline{\hat{\lambda}(\xi)}d\xi.
\]
\begin{comment}
Observe that 
\begin{align*}
&\int \hat{\mu}_\delta(\xi)\overline{\hat{\lambda}(\xi)}d\xi& &=& &\int \hat{\phi}_\delta(\xi)\hat{\mu}(\xi)\overline{\hat{\lambda}(\xi)}d\xi&\\
&& &=& &\int_{B_{K/\delta}(0)}\hat{\phi}_\delta(\xi)\hat{\mu}(\xi)\overline{\hat{\lambda}(\xi)}d\xi+O(\delta^{-n}K^{-N+n})&
\end{align*}
where the $O(.)$ term comes from the fact that $\hat{\phi}_{\delta}$ decays fast off the ball $B_{1/\delta}(0).$
\end{comment}
For each $r>0,$ let $B_r\subset\mathbb{R}^n$ be the metric ball with radius $r$ centred at the origin. For each number $q\geq 0,$ let $S_q=B_{2^{q+1}}\setminus B_{2^{q}}.$ Then we see that,
\begin{align*}\label{I}
 \int |\hat{\phi}_\delta(\xi)\hat{\mu}(\xi)\overline{\hat{\lambda}(\xi)}|d\xi=\int_{B_1}|\hat{\phi}_\delta(\xi)\hat{\mu}(\xi)\overline{\hat{\lambda}(\xi)}|d\xi+\sum_{j\geq 0}\int_{S_j}|\hat{\phi}_\delta(\xi)\hat{\mu}(\xi)\overline{\hat{\lambda}(\xi)}|d\xi.\tag{I}
\end{align*}
Let $\kappa_1<\dim_{l^1}\lambda\leq \dim^I_{l^1}\lambda.$ Observe that
\begin{align*}\label{II}
\int_{S_j}|\hat{\phi}_\delta(\xi)\hat{\mu}(\xi)\overline{\hat{\lambda}(\xi)}|d\xi\ll \frac{1}{2^{j\sigma(M)}}\int_{S_j}|\overline{\hat{\lambda}(\xi)}|d\xi\ll \frac{2^{(j+1)(n-\kappa_1)}}{2^{j\sigma(M)}},\tag{II}
\end{align*}
where $\sigma(M)>0$ can be chosen according to the discussions in Section \ref{sec: fourier manifold}. Thus as long as $\sigma(M)+\dim_{l^1}\lambda>n$, we can choose $\kappa_1$ such that
\[
n-\kappa_1<\sigma(M).
\]
This implies that (we do not take the norm)
\[
\int \hat{\phi}_\delta(\xi)\hat{\mu}(\xi)\overline{\hat{\lambda}(\xi)}d\xi=O(1).
\]
Next, observe that
\[
\lambda(M^{\delta/2})=\int_{M^{\delta/2}} d\lambda(x)\ll \delta^{n-\dim M}\int \mu_\delta(x)d\lambda(x)\ll \delta^{n-\dim M}. 
\]
This finishes the proof of Theorem \ref{Main}.
\subsection{From the counting estimate to an intersection estimate}
Now let $\lambda$ be an AD-regular measure. Suppose that $\dim_{l^1}\lambda$ is large enough so that Theorem \ref{Main} applies.

We then see that
\[
\lambda(M^\delta)\ll \delta^{n-\dim M}.
\]
Since $\lambda$ is AD-regular, we see that $M^\delta\cap \supp(\lambda)$ can be covered by
\[
\ll \delta^{n-\dim M}/\delta^{\Haus \lambda}
\]
many $\delta$-balls. By letting $\delta\to 0$ we see that
\[
\ubox (M\cap \supp(\lambda))\leq \Haus\lambda-(n-\dim M)
\]
if $\Haus \lambda-(n-\dim M)>0.$ Otherwise we have
\[
\ubox (M\cap \supp(\lambda))=0.
\]
This proves Corollary \ref{cor: intersection}.
\subsection{Lower estimate, proof of Theorem \ref{Main2}}\label{sec: lower}
We first reduce the problem a little bit. First, we pick a compact subset of $M$ by performing the smooth truncation discussed in Section \ref{sec: fourier manifold}. This step is not necessary if $M$ is already compact. Let $\mu$ be the chosen smooth compactly supported measure on $M$. Next, we choose an integer $l\geq 0.$ We consider the measure
\[
\mu_l=\frac{1}{p^{nl}}\sum_{d\in \{0,\dots,p^l-1\}^n/p^l} \mu(.+d).
\]
It is an average of $p^{nl}$ many translated copies of $\mu.$ Next, we perform the $\mod \mathbb{Z}^n$ action. We use $\mu_l$ to denote also the image of $\mu_l$ under the action $\mod\mathbb{Z}^n.$ The $\mod \mathbb{Z}^n$ can be performed because $K_{p,D}$ is already $\mathbb{Z}^n$ periodic and $\supp({\mu}_l)$ is compact.\footnote{It is important that $\supp({\mu}_l)$ is compact. Otherwise the image of $\mu_l$ under $\mod\mathbb{Z}^n$ might be supported on a dense subset.} This will let us concentrate on the unit cube $[0,1]^n.$

Now we view the whole situation on $\mathbb{R}^n/\mathbb{Z}^n\approx [0,1)^n.$ It can be checked that for $\xi\in\mathbb{Z}^n$
\[
\hat{\mu}_l(\xi)=\Delta_{p^l|\xi} \hat{\mu}(\xi),
\]
where $\Delta_{p^l|\xi}=1$ if $p^l$ divides all the components of $\xi$ or else $\Delta_{p^l|\xi}=0.$

Let $\delta>0$ be a small number. We consider the thickened measure $\lambda_\delta=\phi_\delta*\lambda_{p,D}.$ Just as $\mu_\delta,$ we see that $\lambda_\delta\asymp \delta^{-(n-\Haus K_{p,D})}$ on $K^{\delta/2}_{p,D}.$ We see that
\begin{align*}\label{III}
\mu_l(K^\delta_{p,D})\gg\delta^{n-\Haus K_{p,D}}\int_{\mathbb{R}^n/\mathbb{Z}^n} \lambda_\delta(x)d\mu_l(x)=\delta^{n-\Haus K_{p,D}}\sum_{\xi\in\mathbb{Z}^n} \hat{\lambda}_\delta(\xi)\hat{\mu}_l(-\xi).\tag{III}
\end{align*}
The above sum converges because $\lambda_\delta$ is a Schwartz function. We can now perform the arguments in Section \ref{sec: L2}. First observe that
\[
\sum_{\xi\in\mathbb{Z}^n} \hat{\lambda}_\delta(\xi)\hat{\mu}_l(-\xi)=\hat{\lambda}_\delta(0)\hat{\mu}_l(0)+\sum_{p^l|\xi,\xi\neq 0} \hat{\lambda}_\delta(\xi)\hat{\mu}_l(-\xi).
\]
For the second sum above, we see that $|\xi|\geq p^l$ because at least one of the components of $\xi$ is non-zero and divisible by $p^l.$ 

We can perform the summation version of the argument in (\ref{I}), (\ref{II}) in Section \ref{sec: L2}. The effect of considering $\mu_l$ instead of $\mu$ is to push the off-zero $L^2$-sum in (\ref{III}) away from the origin.  Since $\dim_{l^1}\lambda_{p,D}\leq \dim^S_{l^1}\lambda_{p,D}.$ We see that
\[
\dim^S_{l^1}\lambda_{p,D}+\sigma(M)\geq \dim_{l^1}\lambda_{p,D}+\sigma(M)>n.
\]
Then we see that (similar to (\ref{II})) as $l\to \infty$
\begin{align*}
& \sum_{p^l|\xi,\xi\neq 0} |\hat{\lambda}_\delta(\xi)\hat{\mu}_l(-\xi)|\\
& \leq \sum_{\xi\in\mathbb{Z}^n, |\xi|\geq p^l} |\hat{\lambda}_\delta(\xi)\hat{\mu}(-\xi)|\\
& \leq\sum_{j\geq l} \sum_{\xi\in\mathbb{Z}^n, p^{j-1}\leq |\xi|< p^j} |\hat{\lambda}_\delta(\xi)\hat{\mu}(-\xi)|\\
& \ll \sum_{j\geq l} p^{-j\sigma(M)} p^{j(n-\kappa_1)},
\end{align*}
where on the last line $\kappa_1$ is a positive number that can be chosen to be arbitrarily close to $\dim_{l^1}\lambda_{p,D}.$ The implicit constant in $\ll$ depends on $\kappa_1$ and does not depend on $l, \delta$.  In particular, $\kappa_1$ can be chosen such that
\[
\kappa_1+\sigma(M)>n.
\]
For this fixed $\kappa_1,$ we see that
\[
\sum_{j\geq l} p^{-j\sigma(M)} p^{j(n-\kappa_1)}\leq p^{l(n-\sigma(M)-\kappa_1)}\frac{1}{1-p^{n-\sigma(M)-\kappa_1}}.
\]
Thus we see that as $l\to\infty,$
\[
\sum_{p^l|\xi,\xi\neq 0} |\hat{\lambda}_\delta(\xi)\hat{\mu}(-\xi)|\to 0
\]
uniformly across $\delta>0.$ Now observe that
\[
\hat{\lambda}_\delta(0)\hat{\mu}_l(0)=1
\]
because $\lambda_\delta, \mu_l$ are all probability measures. Thus as long as $l$ is large enough we have for all small enough $\delta,$
\begin{align*}\label{Lower}
\mu_l(K^\delta_{p,D})>c \delta^{n-\Haus K_{p,D}}\tag{Lower}
\end{align*}
for a constant $c>0$ which depends on the choice of the bump function $\phi.$ 

Observe that $\mu_l$ is essentially a smooth surface measure carried by a manifold. This is not exactly the case because $\mu_l$ is actually carried by a finite union of manifolds. We denote the finite union of manifolds to be $\tilde{M}.$ The estimate
\[
\mu_l(B_\delta(x))\asymp \delta^{\dim M}
\]
holds uniformly for all $x\in \supp(\mu_l)$ and all small enough $\delta>0.$ Of course, the implicit constants depend on $l.$ From here we see that (\ref{Lower}) implies that in order to cover $K^\delta_{p,D}\cap \tilde{M}$ with $\delta$-balls, one needs
\[
\gg \delta^{n-\Haus K_{p,D}-\dim M}=\delta^{-s}
\]
many of them where $s=\Haus K_{p,D}+\dim M-n>0.$

As a simple observation, from (\ref{Lower}), we see that $\tilde{M}\cap K$ cannot be empty. Indeed, if $\tilde{M}\cap K=\emptyset,$ then there is a $\delta_0>0$ such that $d(\tilde{M},K)>\delta_0.$ This is because $M$ and $K$ are closed and compact (we already reduced the whole situation to $[0,1]^n$). This means that $\mu_l(K^\delta_{p,D})=0$ as long as $\delta<\delta_0/10.$ This contradicts (\ref{Lower}). 

Having shown that $\tilde{M}\cap K$ is not empty, we now show that it has positive dimension.  Let $\epsilon>0.$ Let $\delta>0$ be a power of $p^{-1}$ and let $\mathcal{C}$ be any collection of $\delta$-balls with
\[
\#\mathcal{C}\leq (1/\delta)^{\epsilon}.
\]
We want to show that $\mathcal{C}$ cannot cover $\tilde{M}\cap K$ as long as $\epsilon$ is small enough. Let $C\in \mathcal{C}.$ Denote $10C$ by the ball co-centered with $C$ of radius $10\delta.$ Then for each small enough $\delta_C>0,$ there is a number $\eta>0$ such that
\[
\tilde{M}\cap K^{\delta_C}\cap 10C
\]
can be covered with at most 
\begin{align*}\label{eqn: Re-Zoom}
\delta^{\eta}\left(\frac{1}{\delta_C}\right)^{s}\tag{Counting Down}
\end{align*}
many $\delta_C$-balls. This can be seen by rescaling the whole situation by a factor $1/(10\delta).$ Such a zoom action does not change the Fourier decay properties for $K$ or for $M$. Of course, some explicit constants are changed. This results the additional scaling factor $\delta^{\eta}$ for some $\eta>0.$ We will discuss about this later. 

We apply the above argument for each $C\in\mathcal {C}$ and find a small enough number $\delta'>0$ such that each $10C\cap \tilde{M}\cap K$ can be covered by at most
\[
\delta^{\eta}\left(\frac{1}{\delta'}\right)^{s}
\]
many $\delta'$-balls. Thus in total, one needs at most
\[
\delta^{\eta}\left(\frac{1}{\delta}\right)^{\epsilon}\left(\frac{1}{\delta'}\right)^{s}=\delta^{-\epsilon+\eta}\left(\frac{1}{\delta'}\right)^{s}.
\]
many $\delta'$-balls to cover $\bigcup_{C\in\mathcal{C}} 10C\cap \tilde{M}\cap K^{\delta'}.$ If $0<\epsilon<\eta$ we can make $\delta^{-\epsilon+\eta}$ arbitrarily small by choosing $\delta$ to be small. Now we use (\ref{Lower}) for $\delta'.$ We need at least
\[
\gg \left(\frac{1}{\delta'}\right)^{s}
\]
many $\delta'$-balls to cover
\[
\tilde{M}\cap K^{\delta'}.
\]
We thus see that for each small enough $\delta>0,$ as long as $\delta'$ is small enough, there must exist possibly many $\delta'$-balls intersecting $\tilde{M}\cap K$ and they are not intersecting $\cup_{C\in\mathcal{C}} 2C.$ Thus, we apply the above result with $\delta'\to 0$. As a result, we find balls $B_{\delta'}$ of radius $\delta'\to 0$,
\[
B_{\delta'}\cap K\cap \tilde{M}\neq\emptyset
\]
and
\[
d\left(B_{\delta'}, \bigcup_{C\in\mathcal{C}}C\right)>\delta.
\]
Since $K$ and $\tilde{M}$ are compact, we see that there is a point $x\in K\cap \tilde{M}$ with
\[
d\left(x, \bigcup_{C\in\mathcal{C}}C\right)\geq \delta>0.
\]
Thus $\mathcal{C}$ cannot cover $K\cap \tilde{M}.$ We have seen that $K\cap\tilde{M}$ can not be covered by $(1/\delta)^{\epsilon}$ many $\delta$-balls as long as $\delta$ is small enough. This shows that
\[
\lbox K\cap\tilde{M}\geq \epsilon.
\]

Now we discuss how we can choose $\epsilon, \eta$ in (\ref{eqn: Re-Zoom}). Let $\delta>0$ be a power of $p^{-1}.$ We consider a $\delta$-branch of $K,$ say $K'.$ Notice that $K$ is self-similar and the $\delta$-branch $K'$ is just $K$ scaled down by the factor $\delta.$ This scaling procedure affects the Fourier transform as follows. Let $\lambda'$ be the natural missing digits measure supported on $K'$, i.e. $\lambda'=(\lambda(K'))^{-1}\lambda_{|K'}$. Denote $\kappa_2=\Haus K_{p,D}$. Then we see that
\begin{align*}\label{eqn: re-scale}
|\hat{\lambda}'(\xi)|=|\hat{\lambda}(\delta\xi)|.\tag{Rescale}
\end{align*}
We can now argue as in Section \ref{sec: L2} with $\lambda'$ in the place of $\lambda,$
\begin{align*}
&\int_{\mathbb{R}^n} |\hat{\mu}(\xi)\hat{\lambda'}(\xi)|d\xi=\int_{\mathbb{R}^n} |\hat{\mu}(\xi)||\hat{\lambda}(\delta\xi)|d\xi\\
&=\int_{\mathbb{R}^n} |\hat{\mu}(\xi/\delta)||\hat{\lambda}(\xi)|\delta^{-n}d\xi\\
&\ll \delta^{-n}\int_{\mathbb{R}^n} |\xi/\delta|^{-\sigma(M)}|\hat{\lambda}(\xi)|d\xi\\
&=\delta^{-n+\sigma(M)}\int_{\mathbb{R}^n} |\xi|^{-\sigma(M)}|\hat{\lambda}(\xi)|d\xi.
\end{align*}
From here we can use the inequalities in (\ref{II}) to deduce that for $\delta'\to 0,$
\[
\lambda'(M^{\delta'})\ll (1/\delta)^{n-\sigma(M)}{\delta'}^{n-\dim M},
\]
where the multiplicative constant in $\ll$ is the same as at the end of Section \ref{sec: L2}. Since
\[
\lambda(K')\asymp \delta^{\kappa_2},
\]
we see that $K'\cap M^{\delta'}$ can be covered by at most
\[
\ll\frac{(1/\delta)^{n-\sigma}{\delta'}^{n-\dim M}}{\frac{{\delta'}^{\kappa_2}}{\delta^{\kappa_2}}}\ll \delta^{\kappa_2-(n-\sigma)}\left(\frac{1}{\delta'}\right)^{\kappa_2+\dim M-n}
\]
many balls of radius $\delta'.$ Thus we can choose $\eta=\kappa_2-(n-\sigma)>0.$\footnote{Notice that $\kappa_2\geq \dim_{l^1} K$ and by assumption we already have $\dim_{l^1} K+\sigma>n.$} Finally, we can choose any $\epsilon<\eta.$ This shows that
\[
\lbox K\cap\tilde{M}\geq \kappa_2-(n-\sigma).
\]
Finally, if $x\in \tilde{M}\cap K$, then there is a translation vector $d\in \{0,\dots,p^l-1\}^n/p^l$ such that
\[
x+d\in M\cap K_{p,D,l}.
\]
Thus $M\cap K_{p,D,l}\neq\emptyset$ and
\[
\lbox K_{p,D,l}\cap M\geq \kappa_2-(n-\sigma).
\]
To see the latter, observe that
\[
\tilde{M}\cap K\subset K_{p,D,l}\cap \bigcup_{d\in \{0,\dots,p-1\}^n/p^l} (M+d).
\]
For different $d,$ the sets
\[
M_d=K_{p,D,l}\cap (M+d)
\]
are translation copies of each other because $K_{p,D,l}$ is invariant with respect to such translations $d.$ Thus $M_d$ for different $d$ all have the same box covering properties. From here we conclude the lower bound for the lower box dimension of $K_{p,D,l}\cap M.$

For the last statement of Theorem \ref{Main2}, it is enough to observe that for each $\delta>0$ and $\delta$-cube $C$ with $C\cap M\neq\emptyset.$ From the previous argument in this section. There is some possibly large integer $l\geq 0$ depending on $C$ such that
\[
C\cap M\cap K_{p,D,l}\neq\emptyset.
\]

\subsection{A slightly more generalized argument for 'pushing away' the non-zero coefficients}\label{sec: general push}
This section is not needed outside of Section \ref{sec: application}. The reader can skip it for now and come back later. In the previous section, we replaced a measure $\mu$ with
\[
\mu_l=\frac{1}{p^{nl}} \sum_{d\in\{0,\dots,p^l-1\}^n/p^l}\mu(.+d).
\]
The effect of the above averaging process is that for $\xi\in\mathbb{Z}^n,$ $\hat{\mu}_l(\xi)$ is not zero only when the components of $\xi$ are all divisible by $p^l.$ 

Now, we can formulate a slightly more generalized way of performing the above averaging process. Let $p_1,\dots,p_n$ be integers larger than $1.$ Consider the measure
\[
\mu'=\mu_{p_1,\dots,p_n}=\frac{1}{p_1\dots p_{n}}\sum_{d\in \{0,\dots,p_1-1\}p_1^{-1}\times\dots\times\{0,\dots,p_n-1\}p^{-1}_n}\mu(.+d).
\] 
This measure $\mu'$ is an average of $p_1\dots p_n$ translated copies of $\mu.$ The Fourier coefficients $\hat{\mu}'(\xi)$ at $\xi\in\mathbb{Z}^n$ is not zero only when $\xi=(\xi_1,\dots,\xi_n)$ satisfies that for each $i\in\{1,\dots,n\},$
$
p_i|\xi_i.
$
By choosing $p_1,\dots,p_n$ to be all large enough, we again achieve the goal of 'pushing away' the non-zero coefficients as in the previous section.

\subsection{Proof of Theorem \ref{Main3} part 1}
Theorem \ref{Main} has another consequence. Let $M$ be a non-degenerate manifold (or a manifold of finite type) and $K$ be a missing digits set with a large enough $l^1$-dimension such that Theorem \ref{Main} applies. Let $\lambda$ be the corresponding missing digits measure.

From Theorem \ref{Main}, we know that
\[
\limsup_{\delta\to 0} \frac{\lambda(M^\delta)}{\delta^{n-\dim M}}<\infty.
\]
In fact, with extra technical argument, it is likely that the above $\limsup$ can be replaced with $\lim$ but we do not need this. A much easier observation is that one can replace $\limsup$ with $\liminf$ in all the later arguments. We need to use this fact to prove the last assertion of this theorem.

Now we can define a function from $[0,\infty)\times \mathbb{R}^n\times \mathbb{O}(n)$ to $[0,\infty),$
\[
f_{K,M}(t,v,g)=\limsup_{\delta\to 0} \frac{\lambda(T_{t,v,g}(M)^\delta)}{\delta^{n-\dim M}}
\]
where $T_{t,v,g}(M)=t\times g(M)+v,$ i.e., it is the image of $M$ under the rotation $g$, then scaled by $t$ and then translated by $v.$ From here the first part of Theorem \ref{Main3} follows.

Let $\mu$ be a smooth and compactly supported surface measure on $M.$ Then we replace $M^\delta$ with the Schwartz function $\mu_\delta$ (as in Section \ref{sec: L2}) and define the function
\begin{align*}
f_{K,\mu}(t,v,g)=\limsup_{\delta\to 0}\int T_{t,v,g}(\mu_\delta) (x)d\lambda(x).
\end{align*}
For each fixed $\delta>0,$ the quantity in the $\lim$ symbol is continuous viewed as a function with variables $t,v,g.$ Observe that $\mu_\delta$ is a smooth function with $\supp(\mu_\delta)\subset \supp(\mu)^{\delta}$ and  $\mu_\delta\asymp \delta^{-(n-\dim M)}$ on $\supp(\mu)^{\delta/2}.$ From here, we see that
\begin{align*}\label{Asymp}
T_{t,v,g}(\mu_\delta)(x)\asymp \delta^{-(n-\dim M)}\tag{*}
\end{align*}
uniformly for $x\in T_{t,v,g}(\supp(\mu))^{t\delta/2}.$ From here we see that for all $t,v,g$ and all small enough $\delta>0,$
\[
\frac{\lambda((T_{t,v,g}(M))^{t\delta/2})}{\delta^{n-\dim M}}\ll\int T_{t,v,g}(\mu_\delta) (x)d\lambda(x)\ll \frac{\lambda((T_{t,v,g}(M))^{t\delta})}{\delta^{n-\dim M}}.
\]
Thus there is a constant $c>1$ such that
\[
c^{-1} t^{n-\dim M}f_{K,M}(t,v,g)\leq f_{K,\mu}(t,v,g)\leq c t^{n-\dim M} f_{K,M}(t,v,g).
\]

We want to show that $f_{K,\mu}$ is continuous. Then $h=t^{-(n-\dim M)}f_{K,\mu}$ is continuous. In the next section, we will show that $h$ is non-vanishing and this will conclude Theorem \ref{Main3}. 

Observe that
\[
\int T_{t,v,g}(\mu_\delta)(x)d\lambda(x)=\int \hat{T_{t,v,g}(\mu_\delta)}(\xi)\overline{\hat{\lambda}(\xi)}d\xi.
\]
From (\ref{I}), (\ref{II}) we see that the above integrals converge absolutely in a manner that is uniform across $\delta>0$. This is because under the action $T_{t,v,g},$ the Fourier transform is changed accordingly in a manner that preserves the polynomial decay with the same exponent. Thus we can apply (\ref{I}), (\ref{II}). Moreover, if we restrict $t$ within a compact interval away from $0$, then there are constants $c_1,c_2,C_1,C_2>0$ with
\[
C_1\min_{|\xi'|\in [c_1|\xi|,c_2|\xi|]|}|\hat{\mu_\delta}(\xi')|\leq |\hat{T_{t,v,g}(\mu_\delta)}(\xi)|\leq C_2 \max_{|\xi'|\in [c_1|\xi|,c_2|\xi|]}|\hat{\mu_\delta}(\xi')|
\]
for all $\xi.$ The above discussion does not depend on $\delta.$ From here we see that
\[
\limsup_{\delta\to 0}\int T_{t,v,g}(\mu_\delta)d\lambda(x)=\limsup_{\delta\to 0}\int\hat{T_{t,v,g}(\mu_\delta)}(\xi)\overline{\hat{\lambda}(\xi)}d\xi.
\]
The RHS above is continuous viewed as a function with variable $(t,v,g).$ This follows from the following two facts:

a. $T_{t,v,g}$ acts on each of the Fourier coefficients of $\mu_\delta$ continuously in a manner that is uniform across $\delta>0$. Here we emphasize that $T_{t,v,g}$ does not act continuously on the Fourier coefficients in a manner that is uniform across all the frequencies. The continuity here is only pointwise. Thus the role of the cut-off function $\hat{\phi}_\delta$ is important.

b. $\xi\to |\hat{T_{t,v,g}(\mu_\delta)}(\xi)\hat{\lambda}(\xi)|$ is integrable in a manner that is uniform across $\delta>0$ and $t,v,g$ inside any (fixed) compact set $U\subset (0,\infty)\times \mathbb{R}^n\times\mathbb{O}(n).$

From here we deduce that
\[
(t,v,g)\to \limsup_{\delta\to 0}\int T_{t,v,g}(\mu_\delta)(x)d\lambda(x)
\]
is continuous. Thus $f_{K,\mu}$ is continuous. Moreover, the convergence and the continuity is uniform when we restrict $(t,v,g)$ within a compact set $U\subset [0,\infty)\times \mathbb{R}^n\times \mathbb{O}(n)$ away from the set $\{t=0\}\subset [0,\infty)\times \mathbb{R}^n\times \mathbb{O}(n).$ 
\subsection{The integral over group actions: Proof of Theorem \ref{Main3} part 2}\label{sec: group integral}
In order to show that $h=t^{-n+\dim M}f_{K,\mu}$ is non-vanishing, it is enough to show that $f_{K,\mu}$ is non-vanishing (i.e. not identically zero). To show that $f_{K,\mu}$ is non-vanishing, we show that the integral of $f_{K,\mu}$ over a large enough region is positive. Following the previous section, for each $\delta>0,$ consider the integral
\[
\int_U \int T_{t,v,g}(\mu_\delta)(x)d\lambda(x) d(t,v,g),
\]
where $t,v,g$ is integrated over a compact subset $U$ with respect to the Haar measure on $U$. We can exchange the order of the double integral. Observe that for each $x\in\mathbb{R}^n$
\[
\int_U T_{t,v,g}(\mu_\delta)(x)d(t,v,g)\asymp|\{(t,v,g)\in U: T^{-1}_{t,v,g}(x)\in supp (\mu_\delta)\}|/\delta^{n-\dim M}.
\]
where $|.|$ is with respect to the Haar measure restricted on $U$. This measure is not always the Haar measure on $[0,\infty)\times \mathbb{R}^n\times \mathbb{O}(n).$ In fact, it can also be the Haar measure on a subgroup. As $\mu_\delta$ is a compactly supported Schwartz function,
\[
x\to \int_U T_{t,v,g}(\mu_\delta)(x)d(t,v,g)
\]
is non-negative and continuous. Our requirement for $U$ is that $U\cap \{t=0\}=\emptyset$ and for each ball $B\subset\mathbb{R}^n$, each $(t,v,g)\in U$ we have
\begin{align*}\label{Positive}
|\{(t,v,g)\in U: T^{-1}_{t,v,g}(x)\in supp (\mu_\delta)\}|/\delta^{n-\dim M}>\epsilon>0 \tag{Positive}
\end{align*}
for some $\epsilon>0$ and all $x\in B$ in a manner than does not depend on $\delta.$

There are many possible choices of $U.$ 

\subsection*{Case 1:} For example, we can fix $t>0,g\in\mathbb{O}(n)$ and let $v$ range over a sufficiently large ball $B'\subset\mathbb{R}^n$ centred at the origin (i.e. $U=\{t\}\times B'\times \{g\}$). In this way, the Haar measure restricted to $U$ is the Lebesgue measure on the $v$ component restricted to $B'.$ 

\subsection*{Case 2:} Another choice is to fix $v$ and then fix an interval $[a,b]\subset (0,\infty)$ for $t$ with a small enough $a>0$ and a large enough $b>a.$ Finally, we do not restrict $g.$ In this way, the Haar measure restricted to $U$ is the Haar measure of the scaling and rotation group $[0,\infty)\times \mathbb{O}(n).$ In this case, the Haar measure on $U$ is equivalent to the $n$-dimensional Lebesgue measure. 

In both of the cases illustrated above, (\ref{Positive}) is satisfied. In the first example, 
\begin{align*}
|\{(t,v,g)\in U: T^{-1}_{t,v,g}(x)\in supp (\mu_\delta)\}|
\end{align*}
reduces to
\[
|\{v\in B': x+v\in supp (\mu_\delta)\}|=|supp(\mu_\delta)|.
\]
The last term is the Lebesgue measure of $\supp(\mu_\delta)$ which is the $\delta$-neighbourhood of a compact piece of $M.$ It is of order $\delta^{n-\dim M}.$ The second example can be tested via a similar method. 

\subsection*{Case 3:}
We now introduce the third case. We restrict the discussion to $\mathbb{R}^2.$ Let $r>0$ and consider the hyperbola $H_r=\{xy=r\}_{x,y>0}.$ We can choose a piece of $H_r$ by considering $\mu$ to be any compactly support smooth measure on $H_r.$ Let $x\in (0,\infty)^2.$ Consider the line between $x$ and $(0,0)$: $l_x.$ If $l_x\cap \supp(\mu)\neq\emptyset,$ then 
\[
|\{c>0: cx\in supp(\mu_\delta) \}|\asymp \delta,
\]
where $|.|$ is the one dimensional Lebesgue measure. The implicit constants in $\asymp$ depend on $|x|.$ Furthermore, let $r>0$ and $P\subset\mathbb{R}^2$ be compact such that $d(P,(0,0))=r.$ Then the implicit constants in $\asymp$ can be chosen to be the same for all $x\in P.$ 

We have mentioned a few scenarios in which (\ref{Positive}) holds. There are far more situations and we do not provide further examples. Notice that under the Condition (\ref{Positive}),
\[
\int_U f_{K,\mu}(t,v,g)d(t,v,g)=\limsup_{\delta\to 0}\int_U \int T_{t,v,g}(\mu_\delta)(x)d\lambda(x) d(t,v,g)>0.
\]
Here the order of $\lim_{\delta\to 0}$ and the integral $\int_U$ can be changed because that $U$ is compact and $U\cap \{t=0\}=\emptyset$ which implies that the limit is uniform on $U$. Since $f_{k,\mu}$ is continuous, we see that there exists a non-trivial ball $E\subset U$ such that $f_{K,\mu}(t,v,g)>0$ for all $(t,v,g)\in E.$ This concludes the non-vanishing part of Theorem \ref{Main3}.

Finally, for the lower box dimension $\lbox T_{t,v,g}(M)\cap K,$ observe that an estimate like (\ref{Lower}) holds. However, there is a slight difference. We need to replace the above arguments with $\limsup$ being replaced by $\liminf$. Then we have
\[
\lambda((T_{t,v,g}(\supp (\mu)))^\delta)\gg \delta^{n-\dim M}.
\]
This means that $T_{t,v,g}(supp (\mu)^\delta)\cap K$ cannot be covered with $ o(\delta^{n-\dim M-\Haus K})$ many $\delta$-balls. Thus we see that
\[
\mu((T_{t,v,g}(\supp (\mu)))^\delta)\gg \delta^{n-\Haus K}.
\]
From here the rest of the arguments follow as in Section \ref{sec: lower} and we have
\[
\lbox T_{t,v,g}(M)\cap K\geq \Haus K-(n-\sigma(M)).
\]
\begin{comment}
Now, we examine $\mathcal{H}^s(K\cap M\cap B)$ more precisely. Let $\delta>0$ be a small enough number. Let $C_\delta$ be a covering of $K\cap M\cap B$ with balls of radius at most $\delta.$ Without loss of generality, we can assume that $C_\delta$ consists closed cubes with disjoint interiors of side length $p^{-k},k\geq 0, p^{-k}\leq \delta.$ 
Recall that $p$ is the base of $K.$ Moreover, we can assume that for each cube $C$ in $C_\delta,$ either $C\cap K$ only on the boundary of $C$ or else $C\cap K$ is a branch of $K.$ For each $k\geq 0,$ let $C_{\delta,k}$ be the collection of cubes in $C_{\delta}$ of side length $p^{-k}.$ For each cube $C\in C_{\delta,k},$ if $C\cap K$ has positive $\lambda$ measure, then we can apply Theorem \ref{Main} with $\lambda_{C}.$ \end{comment}
\subsection{Sharpness of Theorem \ref{Main} for small missing digits sets: Proof of Theorem \ref{thm: fail}}\label{sec: Fail}
The curve $M$ under consideration is
\[
|x|^k+|y-1|^k=1.
\]
This curve passes through the origin and has an order $k$ contact with the $X$-axis at $(0,0).$ This is the highest order of contact of $M$ with hyperplanes (lines). Loosely speaking, this is because around $(0,0)$ the curve looks like $y=|x|^k$ which has vanishing $(k-1)$th derivative and non-vanishing $k$th derivative at $x=0.$

Let $\delta_0=p^{-l}>0$ for a large integer $l>0.$ Consider the square $[0,\delta_0]^2.$ We decompose this square into smaller squares with side length $p^{-kl}.$ Consider the 'first row' of those smaller squares, i.e. those that intersect the $X$-axis. We see that $M$ intersects each of those smaller squares.
Because the digit set $D$ contains $0.$ This implies that $M^{p^{-kl}}\cap K$ must intersect $\gg (p^{-l}/p^{-kl})^{\Haus K_{p,D_1}}$ many of those squares. 

Let $\delta=p^{-kl}.$ We see that $M^{\delta}\cap K$ cannot be covered with 
\[
\ll \left(\frac{1}{\delta^{(k-1)/k}}\right)^{\Haus K_{p,D_1}}
\]
many $\delta$-balls. Similarly, consider the point $(1,0)$ where $M$ has a $k$th order contact with a vertical line. We see that $M^{\delta}\cap K$ cannot be covered with 
\[
\ll \left(\frac{1}{\delta^{(k-1)/k}}\right)^{\Haus K_{p,D_2}}
\]
many $\delta$-balls. This proves Theorem \ref{thm: fail}. 
\section{Fourier norm dimensions for missing digits measures}\label{sec: l1}
\subsection{General Algorithm}
For the case when $n=1,$ this matter has been discussed in \cite{ACVY}. See also \cite{May2019}. For $n\geq 2$, the arguments change very little. In this paper, we provide some details for being self-contained.

Let $n\geq 1$ be an integer. Let $p>2$ be an integer and $D\subset\{0,\dots,p-1\}^n.$ Consider the missing digit measure $\lambda_{p,D}.$ In this case, $\hat{\lambda}(\xi)$ can be computed with the help of the formula,
\[
\hat{\lambda}_{p,D}(\xi)=\prod_{j\geq 0} \frac{1}{\#D}\sum_{d\in D} e^{-2\pi i (d,\xi)/p^j}.
\]
For convenience, let
\[
g(\xi)=\frac{1}{\#D}\sum_{d\in D} e^{-2\pi i (d,\xi)}.
\]
Then we have
\[
\hat{\lambda}_{p,D}(\xi)=\prod_{j\geq 0} g(\xi/p^j).
\]
We want to estimate for large integers $k\geq 1,$
\[
S_k=\sup_{\theta\in [0,1]^n}\sum_{\xi\in\mathbb{Z}^n,|\xi|_\infty< p^k}|\hat{\lambda}_{p,D}(\xi+\theta)|.
\]
Notice that in the sum, we conditioned on the max norm $|\xi|_\infty=\max\{|\xi_1|,\dots,|\xi_n|\}.$ We now estimate $S_k.$  Let $\theta\in (0,1)^n$ be a vector. Consider the function
\[
f(\theta)=\sum_{\mathbf{i}\in\{0,\dots,p-1\}^n} |g((\mathbf{i}+\theta)/p)|.
\]
Clearly we have for all $\theta,$
\[
0\leq f(\theta)\leq p^n.
\]
Observe that for each $\theta\in [0,1]^n$
\begin{align*}
&S_{k}(\theta)=\sum_{\xi\in\mathbb{Z}^n,|\xi|_\infty< p^{k}}|\hat{\lambda}_{p,D}(\xi+\theta)|\\
&=\sum_{\xi\in\mathbb{Z}^n,|\xi|_\infty< p^{k}}\left|\prod_{j\geq 0}g((\xi+\theta)/p^j)\right|.
\end{align*}
Let $\xi'=\xi+\mathbf{i}$ for some $\mathbf{i}\in p^{k-1}\mathbb{Z}^n.$ Then we have
\[
g((\xi'+\theta)/p^j)=g((\xi+\theta)/p^j)
\]
for all $j=0,1,\dots,k-1.$ From here we see that (recall that $|g|\leq 1$)
\begin{align*}
S_k(\theta)&=\sum_{\xi\in\mathbb{Z}^n,|\xi|_\infty< p^{k}}\left|\prod_{j\geq 0}g((\xi+\theta)/p^j)\right|\\
&\leq \sum_{\xi\in\mathbb{Z}^n,|\xi|_\infty< p^{k}}\left|\prod^{k}_{j= 0}g((\xi+\theta)/p^j)\right|
\\
&= \sum_{\xi\in\mathbb{Z}^n,|\xi|_\infty<p^{k-1}}\sum_{\mathbf{i}\in\{0,\dots,p-1\}^n p^{k-1}} \left|\prod^{k}_{j= 0}g((\xi+\mathbf{i}+\theta)/p^j)\right|
\\
&= \sum_{\xi\in\mathbb{Z}^n,|\xi|_\infty< p^{k-1}}\left|\prod_{j=0}^{k-1} g((\xi+\theta)/p^j)\right|\sum_{\mathbf{i}\in \{0,\dots,p-1\}^n }\left|g(\mathbf{i}p^{-1}+\theta p^{-k}+\xi p^{-k})\right|\\
&\leq \sum_{\xi\in\mathbb{Z}^n,|\xi|_\infty< p^{k-1}}\left|\prod_{j=0}^{k-1} g((\xi+\theta)/p^j)\right|\sup_{\theta'} f(\theta')\\
&\overset{\text{Continue inductively}}{\leq} \sum_{\xi\in\mathbb{Z}^n,|\xi|_\infty< p^{k-2}}\left|\prod_{j=0}^{k-2} g((\xi+\theta)/p^j)\right|(\sup_{\theta'} f(\theta'))^2\\
&\dots\\
&\leq (\sup_{\theta'} f(\theta'))^k.
\end{align*}
Thus we have
\[
S_k(\theta)\leq (\sup_{\theta'}f(\theta'))^{k}.
\]
Therefore we see that
\[
S_k\leq (\sup_{\theta'}f(\theta'))^{k}.
\]
This implies that (we take the Euclidean norm $|\xi|$)
\[
\sup_{\theta\in [0,1]^n}\sum_{\xi\in\mathbb{Z}^n,|\xi|\leq p^k}|\hat{\lambda}_{p,D}(\xi+\theta)|\leq S_k\leq (\sup_{\theta'}f(\theta'))^{k}.
\]
From here one can see that 
\[
n-\frac{\log \sup_{\theta}f(\theta)}{\log p}\leq \dim_{l^1}\lambda_{p,D}.
\]
 \subsection{A crude bound: proof of Theorem \ref{thm: l1 bound} part 1}
 Although it is in principle possible to compute $\dim_{l^1}\lambda_{p,D}$ within any precision, we are not too obsessed with the exact computations. Instead, we provide a rather crude but still useful upper bound for the value
 \[
 \sup_{\theta} f(\theta)
 \]
 when $D$ is a large set. This will give us a lower bound for $\dim_{l^1}\lambda_{p,D}.$
 
 First, observe that
 \[
 \#D g(\xi)=\prod_{j=1}^n\frac{1-e^{2\pi i p\xi_j}}{1-e^{2\pi i \xi_j}}-\sum_{d\notin D}e^{-2\pi i (d,\xi)}.
 \]
 Let $\# D=p^n-t$ for some integer $t>0.$ Then we have
 \[
 -t\leq | (p^n-t) g(\xi)|-\left|\prod_{j=1}^n\frac{1-e^{2\pi i p\xi_j}}{1-e^{2\pi i \xi_j}}\right|\leq t.
 \]
 Now we want to consider the sum
 \[
 f_1(\theta)=\sum_{\mathbf{i}\in\{0,p-1\}^n} |g_1((\mathbf{i}+\theta)/p)|,
 \]
 where
 \[
 g_1(\xi)=\prod_{j=1}^n\frac{1-e^{2\pi i p\xi_j}}{1-e^{2\pi i \xi_j}}.
 \]
 To do this, consider the function $h:\mathbb{R}\to\mathbb{R},$
 \[
 h(x)=\left|\frac{1-e^{2\pi i px}}{1-e^{2\pi i x}}\right|.
 \]
 We want to provide an estimate for
 \[
 H(\theta)=\sum_{j\in\{0,\dots,p-1\}}h((j+\theta)/p).
 \]
 Notice that ($|e^{ix}-1|\leq 2, \forall x\in\mathbb{R}$)
 \begin{align*}\label{eqn: *}
   H(\theta)\leq 2 \sum_{j=0, \{(j+\theta)/p\}\geq 1/p}^{p-1}\frac{1}{|1-e^{2\pi i (j/p)}e^{2\pi i (\theta/p)}|}+\sum_{j: \{(j+\theta)/p\}<1/p} h((j+\theta)/p).\tag{*}
 \end{align*}
 For the first sum in (\ref{eqn: *}), we see that ($|1-e^{2\pi i x}|^2=2(1-\cos(2\pi x))\geq 16x^2$, for $x\in [0,1/2]$)
 \begin{align*}
 &2\sum_{j=0, \{(j+\theta)/p\}\geq 1/p}^{p-1}\frac{1}{|1-e^{2\pi i (j/p)}e^{2\pi i (\theta/p)}|}\leq 8 \sum_{k=1}^{[(p-1)/2]} \frac{1}{|1-e^{2\pi i (k/p)}|}\\
 &\leq \frac{1}{\pi }\sum_{k=1}^{[(p-1)/2]} \frac{k}{p} \leq \frac{p}{\pi} (\log p+1).
 \end{align*}
 For the second sum in (\ref{eqn: *}), notice that there are at most two $j$'s in the sum. As $|h|\leq p$, we have
 \[
 \sum_{j: \{(j+\theta)/p\}<1/p} h((j+\theta)/p)\leq 2p.
 \]
 From here we see that for $p\geq 4,$
 \[
 H(\theta)\leq 2p\log p.
 \]
 We can then use this estimate to see that
 \[
 \sup_{\theta} f_1(\theta)\leq (2p\log p)^n.
 \]
 Thus we have
 \[
 (p^n-t)\sup_{\theta}f(\theta)\leq tp^n+(2p\log p)^n.
 \]
 This implies that
 \begin{align*}\label{Crude bound}
 \dim_{l^1}\lambda_{p,D}\geq n-\frac{\log \frac{tp^n+(2p\log p)^n}{p^n-t} }{\log p}.\tag{Crude bound}
 \end{align*}
 This (\ref{Crude bound}) tells us, for example, if we fix $t>0,$ then as long as $p$ is large enough, $\dim_{l^1}\lambda_{p,D}$ can be arbitrarily close to $n.$ 
 \subsection{Missing digits measures with large bases and rectangular digit sets: proof of Theorem \ref{thm: l1 bound} part 2}
 Let $K_{p,D}$ be a missing digits set. In general, we have
 \[
\frac{1}{2}\Haus K_{p,D} \leq \dim_{l^1} K_{p,D}\leq \Haus K_{p,D}.
 \]
 For missing digits sets, we expect the rightmost side should be closer to the truth. We now make this point clearer.
 
 Let $n\geq 1$ be an integer. Let $p>1$ be an integer. For a missing digits set, we choose a subset $D\subset\{0,\dots,p-1\}^n$ and construct $K_{p,D}.$ We have seen that it is in principle possible to find the value of $\dim_{l^1}\lambda_{p,D}$ with arbitrarily small errors. We can also find a lower bound for $\dim_{l^1}$ with the help of (\ref{Crude bound}).  It turns out that if the digit sets are chosen to be well structured, then we can have a much better estimate than the (\ref{Crude bound}).
 
 Let $D\subset \{0,\dots,p-1\}^n$ be an rectangle, i.e. it is of form
 \[
 [a_1,b_1]\times\dots\times [a_n,b_n]\cap \{0,\dots,p-1\}^n
 \] 
 with integers $a_1\leq b_1,\dots,a_n\leq b_n.$ With this special digit set $D$, we see that the corresponding function $g_{b,D}$ is of form
 	\[
 	g_{b,D}(\theta)=(\#D)^{-1}\sum_{z\in D} e^{-2\pi i((z,\theta))}=(\#D)^{-1}\prod_{j=1}^n e^{-2\pi i a_j \theta_j}\frac{1-e^{-2\pi i (b_j-a_j+1)\theta_j}}{1-e^{-2\pi i\theta_j}}.
 	\]
 	Next, we estimate the sum
 	\[
 	(\#D) f_{b,D}(\theta)=(\#D)\sum_{\mathbf{i}\in\{0,\dots,p-1\}^n}|g_{p,D}(\mathbf{i}+\theta)/p|.
 	\]
 	For each $j\in \{1,\dots,n\},$ define
 	\[
 	S_j(\theta_j)=\left|\frac{1-e^{-2\pi i (b_j-a_j+1)\theta_j}}{1-e^{-2\pi i\theta_j}}\right|.
 	\]
 	Then we see that
 	\begin{align*}
 	(\#D)\sum_{\mathbf{i}\in\{0,\dots,p-1\}^n}|g_{p,D}(\mathbf{i}+\theta)/p|\\
 	=\sum_{\mathbf{i}\in\{0,\dots,p-1\}^n}\prod_{j=1}^n S_j((\mathbf{i}_j+\theta_j)/p)\\
 	=\prod_{j=1}^n \sum_{i\in\{0,\dots,p-1\}} S_j((i+\theta_j)/p).
 	\end{align*}
 	Now we need to estimate for each $j\in\{1,\dots,n\},$
 	\[
 	\sum_{i=0}^{p-1} S_j((i+\theta_j)/p).
 	\]
 	We have already considered this type of sums before, see (\ref{eqn: *}). Notice that $S_j(\theta_j)\leq b_j-a_j+1.$ As a result, we have for $p\geq 4,$
 	\[
 	\sup_{\theta_j} \sum_{i=0}^{p-1} S_j((i+\theta_j)/p)\leq  \frac{p}{\pi}(\log p+1)+2(b_j-a_j+1)\leq 2p\log p.
 	\]
 	Thus we see that
 	\[
 	\sup_{\theta} f_{b,D}(\theta)\leq (\# D)^{-1} (2p\log p)^n.
 	\]
 From here we see that
 \begin{align*}
 \dim_{l^1}\lambda_{p,D}\geq& n-\frac{\log ((\# D)^{-1} (2p\log p)^n)}{\log p}\\
 &= n+\frac{\log \# D}{\log p}-n-\frac{n\log\log p^2}{\log p}\\
 &= \Haus K_{p,D}-\frac{n\log\log p^2}{\log p},
 \end{align*}
 where have used the fact that $\Haus K_{p,D}=\log \#D/\log p.$ 
 
 \section{Proofs of Theorems \ref{thm: special 1}, \ref{thm: special 2}, \ref{thm: special 3}}
 We now prove the special results at the beginning. 
 \subsection{Theorem \ref{thm: special 1}}
 For Theorem \ref{thm: special 1}, one can use Theorems \ref{Main}, \ref{Main2}, \ref{thm: l1 bound}(2) together with the fact that for the Veronese curve in $\mathbb{R}^3,$ the positive number $\sigma$ can be chosen to be $1/3.$ Notice that the missing digits sets in Theorem \ref{thm: special 1} is the threefold Cartesian product $K^3$ where $K$ is the set of numbers on $\mathbb{R}$ whose base $10^{9000}$ expansions only contain digits $\{0,\dots,10^{8100}-1\}.$ With the help of Theorem \ref{thm: l1 bound}(2), it can be checked that $\dim_{l^1}K>8/9.$ Thus we have $\dim_{l^1} K^3>3-3^{-1}.$ The lower bound $1/30$ for box dimension can be derived from the numbers given above.
 \subsection{Theorem \ref{thm: special 2}}
 For Theorem \ref{thm: special 2}, there is still one more step to perform. The missing digits set in consideration is $K\times K$ where $K$ is the set of numbers in $\mathbb{R}$ whose base $10^{9000}$ expansions only have digits in $\{0,\dots,10^{7000}-1\}.$ By Theorem \ref{thm: l1 bound}(2), it can be checked that $\dim_{l^1}$ and then $\dim_{l^1} K\times K> 1.5.$ Let $\lambda$ be the missing digits measures on $K.$
 
 We first illustrate how to show that pinned distance sets have a positive Lebesgue measure. Then we upgrade the result to the non-trivial interval version with the help of Section \ref{sec: group integral}. In this way, we hope to first illustrate clearly the idea behind the proof before the need for too many technical details (zipped in Section \ref{sec: group integral}).
 
 For each $r>0,$ consider the circle $C_r: x^2+y^2=r^2.$ For circles in $\mathbb{R}^2$ (or spheres in $\mathbb{R}^n$), the Fourier decay exponent $\sigma$ can be chosen to be $1/2$ (or $(n-1)/2$ for spheres in $\mathbb{R}^n$).  Theorem \ref{Main} tells us that $\lambda(C^{\delta}_r)\ll_r \delta.$  Moreover, a further insight of Section \ref{sec: L2} tells us that the $\ll_r$ estimate is uniform for $r$ ranging over a bounded interval which is away from zero, i.e. for positive numbers $b>a>0$
 \[
 \lambda(C^\delta_r)\ll_{a,b} \delta 
 \]
 for $r\in [a,b].$ Now we can choose $a>0$ be sufficiently small and $b>0$ being sufficiently large such that 
 \[
 \lambda(B_b(0)\setminus B_{a}(0))>0.\tag{Positive}
 \]
 Consider the set
 \[
 \Sigma_{a,b}=\{r\in [a,b]: C_r\cap K\neq\emptyset\}.
 \]
 Suppose that $\Sigma_{a,b}$ has zero Lebesgue measure. Then we can cover $\Sigma_{a,b}$ with small intervals. More precisely, for each $\epsilon>0,$ there is a countable collection of intervals $I_{j},j\geq 1$ of total length at most $\epsilon$ and they cover $\Sigma_{a,b}.$ Let $I_j$ be one of those intervals. Let $r_j$ be the centre Then we have
 \[
 \lambda(C_{r_j}^{|I_j|/2})\leq c |I_j|,
 \]
 where $c>0$ is a constant depending on $a,b.$ We can sum up the above estimates for all $j,$ 
 \[
 \lambda(B_b(0)\setminus B_{a}(0))\leq c \sum_i |I|_j\leq c\epsilon.
 \]
 As $\epsilon>0$ can be arbitrarily chosen, this contradicts the statement (Positive). Thus we see that
 \[
\Delta_{(0,0)}(K)\supset \Sigma_{a,b}
 \]
 has positive Lebesgue measure. Of course, the choice $(0,0)$ is of no speciality. One can replace it with any other point in $\mathbb{R}^2.$ So that we showed that the pinned distance sets have positive Lebesgue measure.
 
 Now we want to show that the pinned distance sets in fact have non-trivial intervals. This is not very straightforward to show. We use the method introduced in Section \ref{sec: group integral}. Since the circles are compact, we do not need to choose compactly supported smooth surface measures for them. Thus $\mu, \mu_{\delta}$ in Section \ref{sec: group integral} can be simply taken to be the natural ($\delta$-thickened) Lebesgue measures on the circles. For the range of group actions, we take $U=(R^{-1},R)\times \{(0,0)\}\times \mathbb{O}(2)$ for a large enough number $R>0.$ Observe that the circles are invariant under rotations. The arguments in Section \ref{sec: group integral} provide us with a non-vanishing continuous function $f:[R^{-1},R]\to [0,\infty)$ so that whenever $f(r)>0,$ $C_r\cap K\neq\emptyset.$ This shows that $\Delta_{(0,0)}(K)$ contains non-trivial intervals. Again, the point $(0,0)$ is of no significance. One can replace it with any $x\in\mathbb{R}^2.$ However, the value $R$ needs to be changed accordingly. From here the proof of Theorem \ref{thm: special 2} concludes.
 \subsection{Theorem \ref{thm: special 3}}
 Finally, we prove Theorem \ref{thm: special 3}. Consider the class of hyperbola
 \[
 H_r=\{(x,y): xy=r\},r>0.
 \]
 For each $r>0,$ let $(x_1,x_2)\in K\times K\cap (0,1]^2.$ We see that the line connecting $(0,0)$ and $(x_1,x_2)$ will intersect $H_r.$ However, the intersection might be too close to the origin. To overcome this issue, we can consider a branch of $K\times K$ that is away from both of the coordinate lines. Such branches certainly exist, e.g. the image $Y=T_{t,v,g}(K\times K)$ with
 \[
 g=\mathbb{I}, t=10^{-9000}, v=((10^{7000}-1)10^{-9000},(10^{7000}-1)10^{-9000}).
 \]
 Now $Y\subset K\times K$ and $Y$ is far away from the coordinate lines. Let $C_1,C_2>0$ be large enough numbers such that for each $x\in Y,$ the line connecting $(0,0)$ and $x$ intersects $H_r,r\in [C_1^{-1},C_1]$ in $[0,C_2]^2.$
 
 Notice that the $H_r,r>0$ is a class of curves that can be obtained from each via scalings. Now we can apply Section \ref{sec: group integral} (more specifically, the third case) to deduce that 
 \[
 \{r\in [C_1^{-1},C_1]: H_r\cap Y\}
 \]
 contains intervals. From here, the proof of Theorem \ref{thm: special 3} concludes.
 \section{More examples and a question on linear forms}\label{sec: application}
 We explain more applications in this section. For convenience, we fix three missing digits sets on $\mathbb{R}:$
 
 $K_1$: numbers whose base $10^{9000}$ expansions only contain digits in $$\{0,\dots,10^{8100}-1\}.$$ 
 
 $K_2$: numbers whose base $11^{9000}$ expansions only contain digits in $$\{0,\dots,11^{8100}-1\}.$$
 
 $K_3$: numbers whose base $12^{9000}$ expansions only contain digits in $$\{0,\dots,12^{8100}-1\}.$$
 
 The Hausdorff dimensions of $K_1,\dots,K_3$ are equal to $9/10.$ The $l^1$ dimensions of $K_1,\dots,K_3$ are all very close to $9/10$ and in fact they are all larger than $8/9$ by using Theorem \ref{thm: l1 bound}(2).
 
 \begin{exm}
 Consider the hyperbola $\{xy=1\}$ in $\mathbb{R}^2.$ We can apply Theorems \ref{Main}, \ref{Main2} to see that there is an integer $l\geq 0$ and there are infinitely many numbers $t>0$ with
 \[
  10^l t, \frac{10^l}{t}
 \]
 both in $K_1.$
 \end{exm}

We want to list more results that go slightly beyond the scope of Theorem \ref{Main2}. Notice that $K=K_1\times K_2\times K_3$ is not a missing-digit nor a self-similar set by our definition. It is nonetheless AD-regular. Theorem \ref{Main} still applies to this set. However, Theorem \ref{Main2} does not apply to the set $K$. One can review the proof of Theorem \ref{Main2} and find two places where it stops working. 

The first place is at the beginning of Section \ref{sec: lower} where we constructed $\mu_l.$ We need to replace it with the argument in Section \ref{sec: general push}. In our case, we can set $p_1=10^{l_1},p_2=11^{l_2},p_3=12^{l_3}$ for suitable numbers $l_1,l_2,l_3\to\infty$. 

The second place is at the end of Section \ref{sec: lower} where we used (\ref{eqn: re-scale}). It depends on the self-similarity of the underlying set (measure). As our current $K$ is not self-similar, we cannot follow the proof without further modification. 

We now make this modification. Let $\delta>0$ be a small number. We can choose integers $l_1,l_2,l_3$ such that
\[
10^{l_1}\leq \delta^{-1}<10^{l_1+1},
\]
\[
11{l_2}\leq \delta^{-1}<11^{l_2+1},
\]
\[
12^{l_3}\leq \delta^{-1}<12^{l_3+1}.
\]
Then we can decompose $[0,1]^3$ into rectangles of dimension $10^{-l_1}\times 11^{-l_2}\times 12^{-l_3}.$ In this way, we can decompose $K$ into small pieces of the above size. Let $K'$ be one of such pieces. Notice that $K'$ is roughly a box of size $\delta$ up to some bounded multiplicative error. However, $K'$ is not a rescaled copy of $K.$ We can now find the corresponding restricted and normalized measure $\lambda'$ on $K'$ given the original measure $\lambda$ on $K$ which is the product measure of missing-digit measures on $K_1,K_2,K_3.$ Then one can continue the argument in Section \ref{sec: lower}.

After this modification, one can obtain results on products of missing-digit sets (measures). We list three results.

\begin{exm}
	Again consider the hyperbola $\{xy=1\}$ in $\mathbb{R}^2.$ We see that there is an integer $l\geq 0$ and there are infinitely many numbers $t>0$ with $10^l t\in K_1$ and $11^{l}/t\in K_{2}.$
\end{exm}

\begin{exm}
	Consider the Veronese curve $(t,t^2,t^3)_{t\in\mathbb{R}}$ in $\mathbb{R}^3.$ We see that there is an integer $l\geq 0$ and there are infinitely many numbers $t>0$ such that
	\[
	10^t t\in K_1, 11^t t^2\in K_2, 12^t t^3\in K_3. 
	\]
\end{exm}

\begin{exm}
	Consider the curve $\{x^3+y^3=1\}$ on $\mathbb{R}^2.$ For this curve, we can choose $\sigma=1/3.$ There are points $(t_1,t_2)$ on this curve with $10^l t_1\in K_1$ and $11^{l} t_2\in K_2.$
\end{exm}

In those examples, it is possible to study the lower box dimension of the set of points in consideration. This can be done with the same method as in the proof of Theorem \ref{thm: special 1}.

In this paper, we require that the manifold in consideration is of finite type. This excludes the case when it is a line. In fact, despite the simple geometry of lines, their intersections with fractals still remain mysterious:
\begin{ques}
	In Theorems \ref{Main},\ref{Main2} can the manifold $M$ be taken to be irrational hyperplanes? 
\end{ques}
Here, the irrationality condition is crucial. It says that the normal direction of $M$ has rationally independent coordinates. If this condition is not satisfied, then the intersection can be larger than expected, see Section \ref{sec: intersection}. If one allows an $\epsilon$ uncertainty on the exponents then there are satisfactory results. See \cite{Sh}, \cite{Wu}.

\section{Measures with polynomial Fourier decay}\label{sec: fractal measure}
So far, we have only considered problems regarding manifolds intersecting missing digits sets. From the proofs of Theorems \ref{Main}, \ref{Main2}, we see that the property we need for manifolds of finite type is that they support natural surface measures  $\mu$ with whose Fourier transforms have polynomial decay, i.e. for some $\sigma>0$
\[
|\hat{\mu}(\xi)|\ll |\xi|^{-\sigma}.
\]
The proofs of Theorems \ref{Main}, \ref{Main2} can be done with all measures satisfying the above properties. There is no lack of such measures other than surface measures of finite type. We list two special examples. Note that unlike for manifolds with finite type, here the decay exponent $\sigma$ may not be easy to determine and $\mu$ may not be AD-regular. The following example shows that digit expansion and continued fraction expansion are in some sense 'independent'.

\begin{exm}\label{exm: Gibbs}
	Gibbs measures from non-linear interval maps. See \cite{JS}, \cite{CS}. In particular, from \cite[Corollary 1]{JS}, we can deduce the following result:
	
	Let $A\subset\mathbb{N}$ be a finite set. Consider the set $B(A)$ of numbers whose continued fractions only contain digits in $A$. Suppose that $s=\Haus B(A)>1/2.$ Then there is a number $\sigma>0$ such that for missing digits measure $\lambda$ with $\dim_{l^1} \lambda>1-\sigma,$ 
\[
\lambda(B(A)^\delta)\ll \delta^{1-s}.
\]
To prove this result, from \cite{JS} (or \cite{CS}) we know that there is a natural measure $\mu$ supported on $B(A)$ such that $|\hat{\mu}(\xi)|\ll |\xi|^{-\sigma}$ for some $\sigma>0.$ Loosely speaking, this example shows that the continued fraction expansions and digit expansions of numbers are somehow independent in a quantitative way. For example, since $s<1,$ we deduce that $\lambda(B(A)=0.$ This means that $\lambda.a.e$ points $x\in K$ are not in $B(A)$ (This can be also deduced from a much stronger result in \cite{SW19}). 

Next, we want to find points in $B(A)$ that also has missing digits in base $p$ expansion. From the facts that $\dim_{l^1}\lambda>1-\sigma$ and $|\hat{\mu}(\xi)|\ll |\xi|^{-\sigma}$ we deduce that
\[
\int|\hat{\mu}(\xi)\hat{\lambda}(-\xi)|d\xi<\infty.
\]
This implies that $\mu*\lambda$ is absolutely continuous with respect to the Lebesgue measure and moreover the density function is continuous. This implies that the arithmetic sumset $-\supp(\lambda)+\supp(\mu)$ contains intervals. Let $K_{p,D}$ be a missing digits set whose missing digits measure $\lambda$ satisfies $\dim_{l^1}\lambda>1-\sigma.$ Let $I$ be one of those intervals. Let $a\in I$ be any number with terminating base $p$ expansion. Then we see that there are $x\in K_{p,D}$ and $y\in B(A)$ with
\[
-x+y=a.
\]
This implies that
\[
y=a+x.
\]
Since $a$ has terminating base $p$-expansion and $x\in K_{p,D},$ $y$ has missing $p$-adic digits eventually (i.e. there is an integer $l\geq 1$ such that $\{p^la\}\in K_{p,D}$).  
\end{exm}
\begin{exm}
	Patterson-Sullivan measures. See \cite{LNP}. A counting result can be deduced as in Example \ref{exm: Gibbs}. We omit further details.
\end{exm}

 \section{Further discussions}\label{sec: further}
We do not have a complete picture of the distribution of missing digits points around the manifold yet. However, results in this paper provide us with rather satisfactory information in the following situation:
\begin{itemize}
	\item The manifold is sufficiently curved and the missing digits set is sufficiently large with respect to the Fourier $l^1$-dimension.
\end{itemize}

There are several directions for further study:
\begin{itemize}
	
\item First, the largeness of missing digits sets are quantified with the help of $\dim_{l^1}.$ We believe this is not the optimal condition to look at. In fact, we believe that the results in this paper can be proved with $\Haus$ in the place of $\dim_{l^1}.$ For example, we formulate the following conjecture.
\begin{conj}
	Let $n\geq 2$ be an integer. Let $M\subset\mathbb{R}^n$ be a manifold of finite type. Then there is a number $\sigma=\sigma(M)>0$ such that for each missing digits measure $\lambda$ with $\Haus \lambda>\sigma,$
	\[
	\lambda(M^\delta)\ll \delta^{n-\dim M}.
	\]
\end{conj}
Part 2 of Theorem \ref{thm: l1 bound} provides us with examples of missing digits measures that satisfy the conclusion of this conjecture. However, the base of those missing digits measures are all large and the digit sets have to be chosen to be specially structured. Thus the task is to reduce the size of the bases and the structural requirement of the digit sets.

\item Second, what happens if the size of the missing digits set is small? Our theory so far can only be applied when the size of the missing digits set is large.  Then we have obtained an optimal intersection result by combining Theorems \ref{Main}, \ref{Main2}. 

We expect that if the missing digits set has a small enough Hausdorff dimension then it should be rather rare to find those points inside the manifold. We mentioned this point briefly at the beginning of this paper. We formulate the following concrete problem.
\begin{conj}
	Consider the circle $C:x^2+y^2=1.$ For large enough integers $p,$ suppose that $D\subset\{0,\dots,p-1\}^2$ is small enough, say, $\#D\leq 100,$ then $C\cap K_{p,D,l}$ only contains rational points for each $l\geq 0.$
\end{conj}

\item Third, for the method we used in this paper, there are two important factors. We have two sets $A,B$ and we want to study $A\cap B.$ To do this, we find nicely supported measures $\mu_A,\mu_B$ on $A,B$ respectively. Then we need one of the measures, $\mu_A,$ say, to have a power Fourier decay, i.e. for some $\sigma>0,$
\[
|\hat{\mu}(\xi)|\ll |\xi|^{-\sigma}.
\]
For $\mu_B,$ we need that
$\dim_{l^1}\mu_B$ is sufficiently large.

There is no lack of studies of the power Fourier decay property for various measures, e.g. as we have mentioned surface measures carried by manifolds, Gibbs measures, Patterson-Sullivan measures. On the other hand, the study of $\dim_{l^1}$ is relatively new. So far, the best knowledge we have for $\dim_{l^1}\mu_B$ is when $\mu_B$ is a missing digits measure. See also \cite{ACVY} and \cite{Y21}. In particular, in \cite{Y21} a numerical method was proposed to treat self-similar measures which are not missing digits measures. This numerical method does not provide accurate estimates. Apart from these results, almost nothing is known in the case when $\mu_B$ is a general Borel probability measure. We want to ask the following particular problems.
\begin{ques}
Estimate $\dim_{l^1}\mu$ for $\mu$ being:

1. A smooth surface measure carried by non-degenerate manifolds.

2. A self-similar measure with the open set condition.

3. A self-affine measure with the strong irreducibility condition.

4. A self-conformal measure with the open set condition.
\end{ques}
As in \cite{Y21}, answers to this question can help us gain more insights about how rational points are distributed around a specific set, e.g. a self-conformal set, e.g. a Julia set. More generally, as we have discussed in this paper, it is possible to study intersections between different sets from the above list.

\item Theorem \ref{Main2} is not satisfactory because we are only able to find a possibly not sharp lower bound for $\lbox M\cap K_{p,D,l}.$

In fact, as mentioned earlier, we believe that under the hypothesis of Theorem \ref{Main2}, it should be that
\[
\boxd M\cap K_{p,D,l}=\Haus K_{p,D}-(n-\dim M).
\]
We are in fact not too far away from such a result because Theorem \ref{Main2} also tells us that
\[
M^\delta\cap K_{p,D,l}
\]
have the 'correct' size for small enough $\delta>0.$ Thus we see that there are 'enough' points in $K$ which are also close to $M$ but we are not yet able to say that there are 'enough' points in $K_{p,D,l}\cap M.$
\end{itemize}

\section{Acknowledgement}
HY was financially supported by the University of Cambridge and the Corpus Christi College, Cambridge. HY has received funding from the European Research Council (ERC) under the European Union's Horizon 2020 research and innovation programme (grant agreement No. 803711). HY has received funding from the Leverhulme Trust (ECF-2023-186). HY thank P.Varj\'{u} and P. Shmerkin for various comments.

\subsection*{Rights}

For the purpose of open access, the authors have applied a Creative Commons Attribution (CC-BY) licence to any Author Accepted Manuscript version arising from this submission.

\bibliographystyle{amsplain}

\begin{thebibliography}{99}
% 1. Replace 9 by 99 if 10 or more references
% 2. Use "\and" between author names below
\bibitem{ACY} D. Allen, S. Chow and H. Yu, \emph{Dyadic Approximation in the Middle-Third Cantor Set}, preprint: arXiv2005.09300, (2020).

\bibitem{ACVY} D. Allen, S. Chow, P\'{e}ter Varj\'{u} and H. Yu, \emph{Counting rationals and Diophantine approximations in Cantor sets}, work in progress.

\bibitem{B12} V. Beresnevich, \emph{Rational points near manifolds and metric Diophantine approximation}, Annals of Mathematics, \textbf{175}(1), (2012), 187-235.

\bibitem{BDV} V. Beresnevich, H. Dickinson, S. Velani,\emph{Diophantine approximation on planar curves and the distribution of rational points, with an Appendix, Sums of two squares near perfect squares by R.C. Vaughan}. Annals of Mathematics, \textbf{166} (2007), 367-426.



%\bibitem{BDV ref} V. Beresnevich, D. Dickinson and S. Velani, \emph{Measure theoretic laws for lim sup sets}, Mem. Amer. Math. Soc. 179 (2006), no. 846.

\bibitem{BVVZ} V. Beresnevich, R. Vaughan, S. Velani and E. Zorin, \emph{Diophantine Approximation on Manifolds and the Distribution of Rational Points: Contributions to the Convergence Theory}, International Mathematics Research Notices, \textbf{2017}(10), (2017), 2885-2908.


%\bibitem{BV} V. Beresnevich, S. Velani, \emph{A Mass Transference Principle and the Duffin-Schaeffer Conjecture for Hausdorff Measures}, Ann. of Math., \textbf{164}(3), (2006), 971-992.

%\bibitem{B34} A. Besicovitch, \emph{Sets of Fractional Dimensions (IV): On Rational Approximation to Real Numbers}, J. Lond. Math. Soc., \textbf{s1-9}, (1934), 126-131.

\bibitem{B13} J. Bourgain, \emph{Prescribing the binary digits of primes}, Israel Journal of Mathematics, \textbf{194}, (2013),  935-955. 

\bibitem{B15} J. Bourgain, \emph{Prescribing the binary digits of primes  II}, Israel Journal of Mathematics, \textbf{206}(1), (2015),165-182.

\bibitem{Borel} \'{E}. Borel,  \emph{Sur les chiffres d\'{e}cimaux de $\sqrt{2}$ et divers probl\'{e}me de probabilit\'{e}s en cha\^{i}ne}, C. R. Acad. Sci. Paris \textbf{230}, (1950), 591-593. 

\bibitem{BFR11} R. Broderick, L. Fishman and A. Reich, \emph{Intrinsic approximation on Cantor-like sets, a problem of Mahler}, Mosc. J. Comb. Number Theory, \textbf{1}, (2011), 291-300.

\bibitem{BD16}Y. Bugeaud, A. Durand,
\emph{Metric Diophantine approximation on the middle-third Cantor set}, J. Eur. Math. Soc. 18 (2016), no. 6, 1233--1272.

\bibitem{CI} L. De Carli and A. Iosevich, \emph{Some sharp restriction theorems for homogeneous manifolds}, Journal of Fourier Analysis and Applications, \textbf{4}, (1998), 105-128.

\bibitem{CS} S. Connor and T. Sahlsten, \emph{Fourier transform and expanding maps on Cantor sets}, arXiv:2009.01703, (2020).


\bibitem{DMR11} M. Drmota, C. Mauduit and J. Rivat, \emph{The sum-of-digits function of polynomial sequences}, J. Lond. Math. Soc. (2), 84 (1) ,(2011), 81-102. 

\bibitem{FM96} E. Fouvry and C. Mauduit, \emph{Sommes des chiffres et nombres presque premiers}, Mathematische Annalen 305 (3), (1996), 571-600. 

%\bibitem{BD86} J. D. Bovey and M. M. Dodson, \emph{The Hausdorff dimension of systems of linear forms}, Acta
%Arith., \textbf{45}(4), (1986), 337–358.


\bibitem{EFS11} M. Einsiedler, L. Fishman and U. Shapira, \emph{Diophantine approximations on fractals}, Geom. Funct. Anal., \textbf{21}, (2011), 14-35.

%\bibitem{Bern} Bernstein inequality, Encyclopedia of Mathematics.

\bibitem{Fa} K. Falconer, \emph{Fractal geometry: Mathematical foundations and applications, second edition}, John Wiley and Sons, Ltd, (2005).

\bibitem{Fu2}H. Furstenberg, \emph{Intersections of Cantor sets and transversality of semigroups}, Problems in Analysis, Princeton University Press,(1970), 41-59.

%\bibitem{FH09} D. Feng and H. Hu, \emph{Dimension theory of iterated function systems}, Communications on pure and applied Mathematics, \textbf{62}(11), (2009), 1435-1500.
%\bibitem{HF19} R. Fraser and K. Hambrook, \emph{Explicit Salem sets in $\mathbb{R}^n$}, preprint, arxiv:1909.04581, (2020).

%\bibitem{H17} K. Hambrook,  \emph{Explicit Salem sets in $\mathbb{R}^2$}, Adv. Math., 311, (2017), 634–648

%\bibitem{H19} K. Hambrook, \emph{Explicit Salem sets and applications to metrical Diophantine approximation}
%Trans. Amer. Math. Soc., \textbf{371}(6), (2019), 4353–4376.


%\bibitem{JS} T. Jordan and T. Sahlsten, \emph{Fourier transforms of Gibbs measures for the Gauss map},
%Math. Ann., \textbf{364} (3-4), 983-1023, (2016).

%\bibitem{H63} W. Hoeffding, \emph{Probability Inequalities for Sums of Bounded Random Variables},  Journal of the American Statistical Association, \textbf{58}(301), (1963), 13–30.


\bibitem{GIOW}L. Guth, A. Iosevich, Y. Ou and H. Wang, \emph{On Falconer's distance set problem in the plane}, 
Inventiones mathematicae, \textbf{219}, (2020), 779–830.

\bibitem{H14}M. Hochman, \emph{On self-similar sets with overlaps and inverse theorems for entropy}, Annals of Mathematics(2) \textbf{180}(2), (2014), 773-822.

\bibitem{HS12}M. Hochman and P. Shmerkin, \emph{Local entropy averages and projections of fractal measures}, Annals of Mathematics, \textbf{175},(2012), 1001-1059.
	
\bibitem{H20} J. Huang, \emph{The density of rational points near hypersurfaces}, to appear: Duke Mathematical Journal, doi:10.1215/00127094-2020-0004.

\bibitem{HKY} B. Hunt, I. Kan and J. Yorke, \emph{When Cantor sets intersect thickly}, Transactions of the American Mathematical Society \textbf{339}(2), 1993, 869-888.

%\bibitem{HS18} M. Hussain, D. Simmons, \emph{The Hausdorff measure version of Gallagher's theorem – Closing the gap and beyond}, Journal of Number Theory, \textbf{186}, (2018), 211-225.

%\bibitem{H81} J. Hutchinson, \emph{Fractals and Self Similarity}, Indiana Univ. Math. J., \textbf{30}(5), (1981), 713-747.

%\bibitem{K85} J.-P. Kahane, \emph{Some random series of functions}, Cambridge Studies in Advanced
%Mathematics. Cambridge University Press, Cambridge, second edition, 1985.

%\bibitem{K81} R. Kaufman, \emph{On the theorem of Jarn\'{ı}k and Besicovitch}, Acta Arith., \textbf{39}(3), (1981), 265–267.

%\bibitem{K81b} R. Kaufman, \emph{Continued fractions and Fourier transforms}, Mathematika, \textbf{27}, (1981), 262–267.

\bibitem{JS} T. Jordan and T. Sahlsten. \emph{Fourier transforms of Gibbs measures for the Gauss map}. Math. Ann., \textbf{364}(3),(2016),983-1023.

\bibitem{KT} K. Simon and K. Taylor, \emph{Interior of sums of planar sets and curves}, Mathematical Proceedings of the Cambridge Philosophical Society, \textbf{168}(1), (2020), 119-148.


\bibitem{KL20} O. Khalil, M. L\"{u}chi, \emph{Random Walks, Spectral Gaps, and Khinchine's Theorem on Fractals}, arXiv:2101.05797, (2021).

%\bibitem{Khintchine} A. Khintchine, \emph{Einige S\"{a}tze \"{u}ber Kettenbr\"{u}cke, mit Anwendungen auf die Theorie der Diophantischen Approximationen}, Math. Ann. \textbf{92}(1-2), (1924), 115-125.

\bibitem{KLW} D. Kleinbock, E. Lindenstrauss and B. Weiss, \emph{On fractal measures and Diophantine approximation}, Selecta Math. (N.S.) 10 (2004), no. 4, 479--523.

\bibitem{KM98} D. Kleinbock and G. Margulis, \emph{Flows on homogeneous spaces and Diophantine approximation on manifolds}, Ann. of Math., \textbf{148}, (1998), 339-360.

%\bibitem{L98} J. Levesley, \emph{A General Inhomogeneous Jarnik-Besicovitch Theorem}, Journal of Number Theory, \textbf{71}, (1998), 65-80.

\bibitem{LSV ref} J. Levesley, C. Salp and  S. Velani,
\emph{On a problem of K. Mahler: Diophantine approximation and Cantor sets}, Math. Ann. 338 (2007), no. 1, 97--118.

\bibitem{LNP} J. Li, F. Naud, W. Pan,\emph{Kleinian Schottky groups, Patterson-Sullivan measures, and Fourier decay, with an appendix on stationarity of Patterson-Sullivan measures}, Duke Math. J. \textbf{170}(4), (2021), 775-825.

%\bibitem{Mahler84} K. Mahler,
%\emph{Some suggestions for further research}, Bull. Austral. Math. Soc. 29 (1984), no. 1, 101--108.



\bibitem{Ma1} P. Mattila, \emph{Geometry of sets and measures in Euclidean spaces: Fractals and rectifiability}, Cambridge Studies in Advanced Mathematics, Cambridge University Press, (1999).

\bibitem{Ma2} P. Mattila, \emph{Fourier Analysis and Hausdorff Dimension}, Cambridge Studies in Advanced Mathematics, Cambridge University Press, (2015).

\bibitem{MR10} C. Mauduit and J. Rivat, \emph{Sur un probl\'{e}me de Gelfond: la somme des chiffres des
nombres premiers}, Ann. of Math. (2), 171(3), (2010),1591-1646. 

\bibitem{May2019} J. Maynard, \emph{Primes with restricted digits}, Invent. Math. \textbf{217} (2019), 127--218.

\bibitem{MY} C. Moreira and J. Yoccoz, \emph{Stable intersections of regular Cantor sets with large Hausdorff dimensions}, Annals of Mathematics, \textbf{154}(1), (2001), 45-96.

\bibitem{MS18} C. Mosquera and P. Shmerkin, \emph{Self-similar measures: asymptotic bounds for the dimension and Fourier decay for smooth images}, Ann. Acad. Sci. Fenn. Math., \textbf{43}(2) (2018), 823-834.


\bibitem{VV06}B. Vaughan and S. Velani, \emph{Diophantine approximation on planar curves: the convergence theory}, Inventiones mathematicae, \textbf{166}, (2006), 103-124. 

\bibitem{PV2005}
A. Pollington and S. Velani, \emph{Metric Diophantine approximation and `absolutely friendly' measures}, Sel. Math. (N.S.) \textbf{11} (2005), 297--307. 


%\bibitem{G15} S. Gou\"{e}zel, \emph{Limit theorems in dynamical systems using the spectral method}, Proceedings of Symposia in Pure Mathematics, \textbf{89}, (2015), 161-193.



%\bibitem{SS} T. Sahlsten and C. Stevens, \emph{Fourier transform and expanding maps on Cantor sets}, preprint, arxiv:2009.01703, (2020).

\bibitem{Sh} P. Shmerkin, \emph{On Furstenberg's intersection conjecture, self-similar measures, and the $L^q$ norms of convolutions}, Annals of Mathematics, \textbf{189}(2), (2019).

\bibitem{SW19} D. Simmons and B. Weiss, \emph{Random walks on homogeneous spaces and Diophantine approximation on fractals}, Inventiones Mathematicae, \textbf{216}, (2019), 337-394.

\bibitem{Stein} E. Stein, \emph{Harmonic analysis: real-variable methods, orthogonality and oscillatory integrals}, Princeton University Press, 1993.

\bibitem{T} Y. Takahashi, \emph{Products of two Cantor sets}, Nonlinearity, \textbf{30}(5), (2017).

%\bibitem{TV} T. Tao and V. Vu, \emph{Additive Combinatorics}, Cambridge studies in advanced mathematics 105.

%\bibitem{Szusz} P. Sz\"{u}sz, \emph{\"{U}ber die metrische Theorie der Diophantischen Approximation}, Acta. Math. Sci. Hungar,\textbf{9}, (1958),177-193.

%\bibitem{VY20} P. Varj\'{u} and H. Yu, \emph{Fourier decay of self-similar measures and self-similar sets of uniqueness}, preprint, arXiv:2004.09358, (2020).

%\bibitem{Y20} H. Yu, \emph{On the metric theory of inhomogeneous Diophantine approximation: An Erd\H{o}s-Vaaler type result}, preprint, arXiv:2004.05929.

\bibitem{W01} B. Weiss, \emph{Almost no points on a Cantor set are very well approximable}, Proc. R. Soc. Lond. A,\textbf{457}, (2001), 949-952.

\bibitem{Wu}M. Wu, \emph{A proof of Furstenberg's conjecture on the intersections of $\times p$ and $\times q$-invariant sets}, Annals of Mathematics, \textbf{189}(3), (2019), 707--751.



\bibitem{Y20b} H. Yu, \emph{On GILP's group-theoretic approach to Falconer's distance problem}, to appear in Glasgow Math. J., arXiv:1810.00987

\bibitem{Y21} H. Yu, \emph{Rational points near self-similar sets}, preprint,  arXiv:2101.05910, (2021).


\bibitem{ZRZJ} B. Zhao, X. Ren, J. Zhu and K. Jiang, \emph{Arithmetic on self-similar sets}, Indagationes Mathematicae
\textbf{31}(4), (2020), 595-606.

\end{thebibliography}

\end{document}